\documentclass[11pt,sort]{elsarticle}

\usepackage[bottom = 2in]{geometry}

\usepackage{amsmath,amssymb,amsthm,mathtools}
\usepackage{enumerate}
\usepackage{parskip}
\usepackage{cancel}
\usepackage{subfigure}

\newcommand{\isDraft}{false}

\newcommand{\mcf}{\mathcal}
\newcommand{\mbf}{\mathbf} 
\newcommand{\mbb}{\mathbb}
\renewcommand{\Re}{\mbb{R}}
\newcommand{\eqdef}{:=}

\newcommand{\tpose}{^{\smash{\sf{T}}}} 
\newcommand{\set}[2]{\left\{ #1\ \left| \ #2 \right. \right\}}
\newcommand{\norm}[1]{\left\lVert#1\right\rVert}
\newcommand{\Reynolds}{\ensuremath{Re}}

\newcommand{\ReynoldsEnergyLim}{\Reynolds_{\rm{e}}}

\newcommand{\nnz}{\operatorname{nnz}}

\newcommand{\innerprod}[2]{\left\langle{#1},{#2}\right\rangle}
\newcommand{\half}{\frac{1}{2}}
\newcommand{\nablaa}{\nabla_a}

\def\hsmash{\mathpalette\mathrlapinternal} 
\def\mathrlapinternal#1#2{%
\rlap{$\mathsurround=0pt#1{#2}$}}

\begin{document}

\title{Global Stability Analysis of Fluid Flows using Sum-of-Squares}
\author{P.~Goulart\corref{cor1}}
\address{Automatic Control Laboratory, ETH Z\"urich,\\ Physikstrasse 3, 8092 Zurich, Switzerland }
\author{S.~Chernyshenko\corref{cor2}}
\address{Department of Aeronautics, Imperial College London, \mbox{Prince Consort Road, London, SW7 2AZ, United Kingdom}}
\cortext[cor1]{\texttt{pgoulart@control.ee.ethz.ch}~(Corresponding author). }
\cortext[cor2]{\texttt{s.chernyshenko@imperial.ac.uk}}


\begin{abstract}
This paper introduces a new method for proving global stability of fluid flows through the construction of Lyapunov functionals.  For finite dimensional approximations of fluid systems, we show how one can exploit recently developed optimization methods based on sum-of-squares decomposition to construct a polynomial Lyapunov function.  We then show how these methods can be extended to infinite dimensional Navier-Stokes systems using robust optimization techniques.  Crucially, this extension requires only the solution of infinite-dimensional linear eigenvalue problems and finite-dimensional sum-of-squares optimization problems.  

We further show that subject to minor technical constraints, a general polynomial Lyapunov function is always guaranteed to provide better results than the classical energy methods in determining a lower-bound on the maximum Reynolds number for which a  flow is globally stable, if the flow does remain globally stable for Reynolds numbers at least slightly beyond the energy stability limit.  Such polynomial functions can be searched for efficiently using the SOS technique we propose.

\end{abstract} 

\maketitle

\section{Background and problem statement}\label{sec:intro}

In this paper we propose a new analytical method for determining whether a fluid flow is globally stable.  This new approach has its origins in two hitherto distinct research areas. The first of these is the classical energy approach of \cite{Serrin59,Joseph76}, which provides conservative lower bounds on the stability limits of flows by analyzing the time evolution of the energy of flow perturbations.  The other is the emerging field of sum-of-squares (SOS) optimization over polynomials, which can be used to prove global stability of finite-dimensional systems of ordinary differential equations with polynomial right-hand sides~\cite{Parrilo00, Parrilo2003}.   It is our hope that the present text is written in such a way that it will be understandable to researchers from either of these two areas, which until now have remained almost completely isolated. We are aware of only two other publications where SOS methods have been used to analyze the behavior of similar systems governed by partial differential equations, namely \cite{PP06,HKI08}.

We apply SOS methods to the problem of assessing global stability of an incompressible flow. The velocity  $\mbf w$ and pressure $p'$ of a flow of viscous incompressible fluid, evolving inside a bounded domain $\Omega$ with boundary $\partial \Omega$ under the action of body force $\mbf{f}$, is governed by the Navier-Stokes and continuity equations
\begin{subequations}\label{eqn:NS}
\begin{align}
\frac{\partial \mbf{w}}{\partial t} + \mbf{w}\cdot\nabla\mbf{w}
&= -\nabla p' + \frac{1}{\Reynolds}\nabla^2\mbf{w} + \mbf{f} \label{eqn:NS:a}\\
\nabla \cdot \mbf{w} &= 0, \label{eqn:NS:b}
\end{align}
\end{subequations}
with a boundary condition $\mbf{w}=0$ on $\partial\Omega.$ Here $\Reynolds$ is the Reynolds number, which is a dimensionless parameter indicating the relative influence of viscous and inertial forces in the flow.  In what follows we will make extensive use of an inner product of vector fields defined as
$$
\innerprod{\mbf{u}}{\mbf{v}} \eqdef \int_\Omega \mbf{u}\cdot\mbf{v}\,d\Omega
$$
with the usual $\mcf{L}_2$ norm $\norm .$ defined as
$
\norm{\mbf{u}}^2 \eqdef \innerprod{\mbf{u}}{\mbf{u}}$.  Similarly, we define (using standard Einstein summation notation throughout)
\[
\innerprod{\nabla \mbf{u}}{\nabla \mbf{v}} \eqdef \int_\Omega \frac{\partial u^i}{\partial x^j} \frac{\partial v^i}{\partial x^j} d\Omega 
\]
 and $\norm{\nabla \mbf{u}}^2 \eqdef \innerprod{\nabla \mbf{u}}{\nabla \mbf{u}}$.

We say that a steady solution $\mbf{w}=\bar{\mbf u},$ $p'=\bar p$ of the system \eqref{eqn:NS} is globally stable if, for each $\epsilon > 0$, there exists some $\delta > 0$ such that $\norm{\mbf w - \bar{\mbf u}} \le \delta$ at time $t_0$ implies that $\norm{\mbf w - \bar{\mbf u}} \le \epsilon$ for all time $t \ge t_0$.   We say that it is globally asymptotically stable if in addition
$\mbf{w}\to\bar{\mbf u}$ as time $t\to\infty$ for any initial conditions.  For the system \eqref{eqn:NS}, these stability conditions ensure laminar flow.   Our principal aim is to identify the largest value~$\Reynolds$ for which these conditions can be guaranteed to hold for the system \eqref{eqn:NS}.

The results described in this paper can be extended to other types of boundary conditions, most notably to the frequently encountered case of periodic boundary conditions.  A useful property of systems with such boundary conditions is that, for any solenoidal vector field or fields satisfying these boundary conditions and for any scalar function $\phi$, the following hold true:
 $\innerprod{\mbf{v}}{\nabla\phi}=0,$
$
\innerprod{\mbf v_1}{\nabla^2\mbf v_2} =\innerprod{\mbf v_2}{\nabla^2\mbf v_1}= -\innerprod{\nabla \mbf v_1}{\nabla \mbf v_2}$
and
 \begin{equation}\label{eqn:tripleproduct}
\innerprod{\mbf{v}_1}{\mbf{v}_2\cdot\nabla\mbf{v}_3}=-\innerprod{\mbf{v}_3}{\mbf{v}_2\cdot\nabla\mbf{v}_1},
\end{equation} 
hence
 \begin{equation}
\innerprod{\mbf v}{\mbf v \cdot\nabla\mbf v}=0.\label{eqn:EnergyConservation}
\end{equation}
These properties can be proved by applying standard identities from vector calculus or simply by integrating by parts, using incompressibility ($\nabla \cdot \mbf v=0$) and applying the boundary conditions. These properties are, of course, well known.
 
In particular, the identity \eqref{eqn:EnergyConservation} plays an important role in the energy approach to proving global stability \cite{Serrin59}. 
Defining velocity  perturbations $\mbf{u}\eqdef\mbf{w}-\bar{\mbf u}$ and pressure perturbations $p\eqdef p'-\bar p$, the system
\eqref{eqn:NS} can be written as
\begin{subequations}\label{eqn:NSpert}
\begin{align}
\frac{\partial \mbf{u}}{\partial t} + \mbf{u}\cdot\nabla\mbf{u} + S(\mbf u,\bar{\mbf u})
&= -\nabla p + \frac{1}{\Reynolds}\nabla^2\mbf{u}
\label{eqn:NSpert:S}\\
\nabla \cdot \mbf{u} &= 0 \label{eqn:NSpert:b},
\end{align}
\end{subequations}
where  $S(\mbf u, \mbf v) \eqdef \mbf u \cdot \nabla \mbf v + \mbf v \cdot \nabla \mbf u$ is introduced for compactness of notation. 
The time rate of change of the velocity perturbation energy can be obtained by taking the inner product of both sides of   
\eqref{eqn:NSpert:S} with $\mbf{u}$, to obtain the energy equation
\begin{equation}\label{eqn:EnergyEquation}
 \frac{\partial\norm{\mbf u}^2/2}{\partial t}=\frac{1}{\Reynolds}\innerprod{\mbf u}{\nabla^2\mbf u}
- \innerprod{\mbf u}{S(\mbf u,\bar{\mbf u})}.
\end{equation}
Note that the nonlinear term $\mbf u \cdot \nabla \mbf u$ in \eqref{eqn:NSpert:S} does not feature in the energy equation because of the identity \eqref{eqn:EnergyConservation}. This is particularly useful because it allows one to obtain immediately an (albeit conservative) method for checking stability of the system \eqref{eqn:NSpert}, which we now describe~briefly.

There exists a real constant $\kappa$ such that for all solenoidal $\mbf u$, 
\begin{equation}\label{eqn:EnergyInequality}
\frac{1}{\Reynolds}\innerprod{\mbf u}{\nabla^2\mbf u}
- \innerprod{\mbf u}{S(\mbf u,\bar{\mbf u})}  \le \kappa \norm{\mbf{u}}^2.
\end{equation}

We can follow the procedure in~\mbox{\cite[p.\ 33-34]{DB95}} to find the smallest such $\kappa$ via solution of an eigenvalue problem.
The smallest $\kappa$ satisfying \eqref{eqn:EnergyInequality} is the solution to the optimization problem
\[
\sup 
\left(\frac{-
\frac{1}{\Reynolds}\norm{\nabla \mbf u}^2 - \int_\Omega \mbf{u}\cdot(\nabla \bar{\mbf u})\cdot \mbf {u}\,dV
}
{\norm{\mbf u}^2}
\right),
\]
subject to the incompressibility condition $\nabla \cdot \mbf u = 0 $ and the boundary conditions. Since the objective function in this optimization problem is homogeneous in $\mbf u$, one is free to optimize over the numerator only, with additional constraint $\norm{ \mbf u} = 1$.  This leads to the eigenvalue problem
\begin{equation}\label{eqn:EigenEnergy}
\begin{gathered}
\lambda \mbf u   = (e - \frac{1}{\Reynolds}\nabla^2)\mbf u + 
\nabla p \\
\nabla \cdot \mbf u = 0, \,\, \left.\mbf u\right|_{\partial\Omega}=0,
\end{gathered}
\end{equation}
where $p$ is the Lagrange multiplier for the incompressibility condition, $\lambda$ is the Lagrange multiplier for the unit norm condition $\norm{\mbf u}=1$ and $e$ is the base-flow rate of strain tensor with components
\[
e_{ij}(x) \eqdef \frac{1}{2}\left[\frac{\partial \bar u^i}{\partial x^j}
+ \frac{\partial \bar u^j}{\partial x^i}\right ].
\]
Note the identities used in arriving at \eqref{eqn:EigenEnergy}; $
\innerprod{\mbf u}{S(\mbf u,\bar{\mbf u})} = 
\cancel{\innerprod{\mbf u}{\bar{\mbf u} \cdot \nabla \mbf u}}^{{\,}^{=0}}\!\!
+\innerprod{\mbf u}{\mbf u \cdot \nabla \bar{\mbf u}}$ and $\mbf u\cdot(\nabla \mbf{\bar u})\cdot\mbf u=\mbf u\cdot e\cdot\mbf u.$ Since the operator in \eqref{eqn:EigenEnergy} is self-adjoint \cite{DB95}, all the eigenvalues $\lambda_k$ of \eqref{eqn:EigenEnergy} are real.  If these eigenvalues are ordered by decreasing value with $\lambda_1$ being the largest, then the inequality \eqref{eqn:EnergyInequality} is tight with $\kappa = \lambda_1$.
 If the largest eigenvalue $\lambda_1<0$, then \eqref{eqn:EnergyEquation} is always negative for $\norm{\mbf u} \neq 0$ and hence the energy $\norm{\mbf u}^2/2$ is a Lyapunov functional for \eqref{eqn:NSpert}, thus proving the global stability.

A particularly nice feature of the energy approach is that proving global stability by this method requires only the solution of  linear eigenvalue problems, even though \eqref{eqn:NSpert} is a system of nonlinear partial differential equations.  This is a direct consequence of the unique advantages of using $\norm{\mbf u}$ as a Lyapunov functional, in particular the opportunity to exploit~\eqref{eqn:EnergyConservation}.  In general, any other choice of Lyapunov functional would result in a stability condition featuring the nonlinear inviscid term $\mbf u \cdot \nabla \mbf u$, in contrast to the more benign energy stability condition \eqref{eqn:EnergyEquation}.

On the other hand, 
the energy approach can give very conservative results, in the sense that the largest $\Reynolds$ for which global stability can be proven by this method is generally well below the maximum $\Reynolds$ for which the flow is generally observed to be globally stable, either numerically or experimentally.  
 
The approach proposed in the present study aims to improve this bound, using a partial Galerkin decomposition of the infinite dimensional system~\eqref{eqn:NSpert}.  Finite dimensional methods, based on recently developed techniques in polynomial optimization, are used to define a Lyapunov functional that is nonlinear in a finite number of terms, while otherwise maintaining some of the attractive numerical advantages of energy methods for the remaining (infinite dimensional) dynamics.  We stress that our results suggest a way of computing a  Lyapunov functional verifying stability of the infinite dimensional system \eqref{eqn:NSpert}, and not some truncated finite dimensional approximation thereof. 

\subsection{Finite Dimensional Systems and the Sum-of-Squares Decomposition}

We first comment briefly on the state of the art in direct methods for computing Lyapunov functions for finite dimensional nonlinear systems.    Suppose that the evolution of a finite-dimensional system with state vector $a \in \Re^n$ is governed by a set of ordinary differential equations (ODEs)
\begin{equation}\label{eqn:genericODE}
\dot a=f(a)
\end{equation}
with equilibrium point $a=0$.   We will use $\nablaa$ throughout to indicate the gradient of a scalar function defined on this $n$-dimensional state space, and otherwise use $\nabla$ to indicate the gradient or divergence of functions in physical space, as in \eqref{eqn:NS}.

The origin of the system \eqref{eqn:genericODE} is globally asymptotically stable if there exists a continuously differentiable \emph {Lyapunov function} $V : \Re^n \to \Re$ such that $V(0) = 0$, $V(a) > 0$ for all \mbox{$a \in \Re^n \setminus \{0\}$} and $\dot V(a) = \nablaa V \cdot f(a)  < 0$  for all $a \in \Re^n \setminus \{0\}$  \cite{Khal02}.   Given a Lyapunov candidate function $V(a)$ and associated $\dot V(a)$, these conditions amount to checking global positivity or negativity of functions.   There is no general method for performing such a check, nor any systematic way of constructing Lyapunov functions for general systems of ODEs.  

A truncated Galerkin approximation reduces the Navier-Stokes equations \eqref{eqn:NSpert} to a system of ODEs in exactly the form \eqref{eqn:genericODE}, but with polynomial (in fact quadratic) right hand side~$f(a)$.  In this particular case, checking that a polynomial function $V(a)$ serves as a Lyapunov function reduces to verifying the positive-definiteness of the two related polynomials $V(a)$ and $-\dot V(a)$.  However, verifying positive-definiteness of a general multivariate polynomial is still NP-hard in general, and is a classical problem in algebraic geometry.

Nevertheless, there has been significant recent progress in stability analysis of polynomial systems using \emph{sum-of-squares} optimization methods, which were first employed in the context of dynamical systems in \cite{Parrilo00}.   These methods are based on a recognition that  a {sufficient} condition for a polynomial function to be positive-definite is that it can be rewritten as a sum-of-squares (SOS) of lower order polynomial functions\footnote{This condition is not a necessary one however, apart from certain exceptional cases involving relatively few variables or low-order polynomials.  An example of a positive-definite polynomial function that is \emph{not} a sum-of-squares is the Motzkin polynomial~\mbox{$y^4z^2 + y^2z^4 - 3y^2z^2 + 1$}.  A nice account of the history of this problem, the 17th of Hilbert's 23 famous problems posed at the turn of the 20th century, can be found in~\cite{Reznick2000}.}.  
Verifying this stronger condition, and solving other problems related to such representations, is significantly simpler than verifying global positivity in general. Therefore the general approach of sum-of-squares optimization in control applications is to search for a Lyapunov function $V(a)$ and associated function $\dot V(a)$ that satisfy sum-of-squares conditions.

Every polynomial $V(a)$ of order $2k$ can be represented as a quadratic form of monomials of order less than or equal to $k$, i.e.\ in the form $V=C_{ij}m_i(a)m_j(a)$. The monomials $m_i(a)$ in this factorization are expressions of the form $a_1^{k_1}a_2^{k_2}\dots a_n^{k_n}$, with integer exponents $k_i\ge 0$ satisfying $\sum_i^nk_i\le k$. If the matrix $C_{ij}$ is positive-definite and symmetric (the latter is always possible), then it can be diagonalized by a suitable linear transformation of the monomial set. Since all the diagonal elements in the resulting expression will be positive, this gives a representation of the polynomial as a sum of squares of polynomials of lower order.  Such a representation is known as a sum-of-squares decomposition.

Hence, the problem of finding a Lyapunov function is reduced to finding coefficients $C_{ij}=C_{ji}$ such that this matrix is positive-definite and the corresponding matrix factorization representing $-\dot V(a)$ is also positive-definite.  The relationship between the coefficients of $V(a)$ and the matrix $C_{ij}$ amounts to a set of linear equality constraints, with similar linear equality constraints relating the coefficients of $\dot V(a)$ and its factorization. A further set of equality constraints couple the coefficients of the polynomial functions $V(a)$ and~$\dot V(a)$.  

Problems such as that described above can be solved efficiently since the set of positive-definite matrices is convex.  
The general field of optimization theory and numerical methods related to such problems is known as \emph{semidefinite programming}%
\footnote{
Strictly speaking, the problem described here has positive-definite matrix constraints, rather than semidefinite constraints as in standard semidefinite programming.  The conditions described above can be recast as semidefinite constraints via inclusion of appropriate terms, e.g. $V(a) > 0$ if $V(a) - \epsilon\norm{a}^2 = C_{ij}m_i(a)m_j(a) \ge 0$ for some small $\epsilon > 0$.  In this case the semidefiniteness condition $C_{ij} \succeq 0$ is sufficient.}, 
and such problems are solvable in an amount of time that is polynomial in the size of their problem data \cite{Wolkowicz00,Todd01,Boyd04}.  Standard software tools are freely available  for posing and solving sum-of-squares problems \cite{yalmip,YalmipSOS09,PP05} as semidefinite programs.  

\subsection{Application of SOS methods to fluid systems}

With respect to the SOS approach, ODE systems obtained via finite-dimensional approximation of the Navier-Stokes equations require special treatment for two reasons. First, since the ODEs describing the dynamics are quadratic, if $V(a)$ is of even degree (as it must to be positive-definite), then $-\dot V(a)$ is formally of odd degree, and hence will not be positive-definite in general. 


The second reason is more subtle. Consider the behaviour of the Lyapunov functional for very large values of $\norm{\mbf u}$. In this case the term $S(\mbf u,\bar{\mbf u})$ and the viscous term $\frac{1}{\Reynolds}\nabla^2\mbf{u}$ in \eqref{eqn:NSpert:S} become small relative to the nonlinear term $\mbf{u}\cdot\nabla\mbf{u}$, and the dynamics become approximately inviscid. In the limit the time derivative of a Lyapunov functional can remain negative or can tend to zero. In the first case the high-$\norm{\mbf u}$ asymptotics of the Lyapunov functional would be a Lyapunov functional of the zero solution for the inviscid flow with zero forcing. This is impossible because in such an inviscid flow the energy is conserved, i.e.\ the flow does not decay to rest. In the second case the asymptotics will be a functional remaining constant on any solution, i.e.\ it will be an invariant of the inviscid flow. 

In the 2D case there are an infinite number of such invariants, so that no conclusions can be drawn, but in the 3D case the only invariant known so far is energy.  Hence, it is highly likely that the high-$\norm{\mbf u}$ asymptotics of the Lyapunov functional of the viscous flow match $\norm{\mbf{u}}$ itself or a monotone function thereof. Therefore, it is reasonable to limit the search for a Lyapunov functional to functionals with such asymptotics. Moreover, since $\norm{\mbf w}$ has the same high-$\norm {\mbf u}$ asymptotics as $\norm{\mbf u}$ itself and since, unlike $\norm {\mbf u},$ $\norm {\mbf w}$ decays monotonously in a viscous flow for any $Re$ if $\norm {\mbf u}$ is large enough, it is sensible to seek a Lyapunov functional that behaves like $\norm {\mbf w}$ at large $\norm {\mbf u}$. Although this argument is not rigorous, it helps in guessing the structure of an appropriate Lyapunov functional. A more technical analysis of the asymptotic behaviour of Lyapunov functionals is given in Section~\ref{sec:extremeNorm_a} for the finite-dimensional case.  

We describe in Sections \ref{sec:finiteDimProblems} and \ref{sec:FiniteDim} how one can apply SOS methods to truncated ODE approximations to \eqref{eqn:NSpert} to obtain numerical estimates of the maximum value $\Reynolds$ for which the system is global stable.  We supply a numerical example illustrating the application of these methods in Section \ref{sec:example}, where a stability limit approximately seven times larger than the value demonstrable via energy methods is obtained, and which is close to the global stability limit estimated by direct numerical simulation.    A preliminary version of these results was also presented in~\cite{GC10}.

There remains the question of convergence of global stability results obtained in this way as the number of modes in the truncated Galerkin approximation tends to infinity; global stability of a truncated approximation does not imply global stability of the Navier-Stokes solution.   This problem is particularly acute since one cannot realistically expect to apply SOS methods to high resolution approximations of the Navier-Stokes equations, since the size of the related optimization problems quickly becomes unmanageable.  

We therefore demonstrate in Sections~\ref{sec:infDimProblems} and~\ref{sec:InfiniteDim} how these difficulties can be overcome by searching for a Lyapunov functional of the Navier-Stokes system in the form $V=V(a,q^2),$ where $a$ is a finite-dimensional vector of the amplitudes of several Galerkin modes, and $q^2$ is the collective energy of all of the remaining modes. Such an approach requires estimates via $q$ for the terms stemming from the nonlinearity of the Navier-Stokes system, which is immediately reminiscent of the difficulties preventing proof of existence of solutions to the Navier-Stokes equations. 
In the present context it turns out, however, that the required estimates are available since they are needed only for the effect of higher-order modes on the finite set of Galerkin modes. As a result, the required estimates can be obtained by solving only linear eigenvalue problems in infinite dimensions and certain maximization problems in finite dimensions, and the resulting system can be treated using the SOS approach. 

We show in Section \ref{sec:energyCompare} that, with a suitable basis for the Galerkin approximation, the proposed approach is guaranteed always to give results at least as good as the standard energy approach.  We further show in \ref{app:existence} that if the  flow remains globally stable in some range of $\Reynolds$ beyond the maximum $\Reynolds$ for which stability can be proved using the energy approach, then a polynomial Lyapunov function is still guaranteed to exist in at least part of this extended range.

\section{Finite Dimensional Flow Models}\label{sec:finiteDimProblems}

We assume throughout that the perturbation velocity $\mbf{u}$ can be written as 
\begin{equation}\label{eqn:ansatz}
\mbf {u}(x,t) = a_i(t)\mbf u_i(x) + \mbf{u}_s(x,t), \qquad i=1,\dots,n
\end{equation}
where the basis functions $\mbf u_i$ are  mutually orthogonal, solenoidal and satisfy the boundary conditions.  Likewise, $\mbf{u}_s$ is assumed to be solenoidal, to satisfy the boundary conditions,  and to be orthogonal to the bases $\mbf{u}_i$.   We assume also that each of the basis functions has unit norm, i.e.\ $\norm{\mbf u_i} = 1$. For brevity, we denote by $\mcf S$ the set of all possible vector fields that are solenoidal, satisfy the boundary conditions and are orthogonal to all $\mbf{u}_i$, so that $\mbf u_s \in \mcf S$. 

In order to address the global stability of the nonlinear Navier-Stokes system \eqref{eqn:NSpert}, we will partition its dynamics into the interaction of an ODE, representing the evolution of the basis weights $a$, and a PDE, representing the remaining unmodeled modes of the system~$\mbf u_s$.   We work initially with the ODE part only, and hence assume initially that~$\mbf{u}_s= 0$.   

First substitute \eqref{eqn:ansatz} into \eqref{eqn:NSpert:S} and take an inner product of both sides with each of the basis functions $\mbf{u}_i$ in turn, yielding an ODE in the form
\begin{equation}\label{eqn:nominalODE}
\dot a_i + \innerprod{\mbf u_i}{\mbf u_j \cdot \nabla \mbf u_k}a_ja_k + \innerprod{\mbf u_i}{S(\mbf u_j,\bar{\mbf u})}a_j
= \frac{1}{\Reynolds}\innerprod{\mbf u_i}{\nabla^2\mbf{u}_j}a_j.
\end{equation}
Defining matrices $\Lambda$, $W$ and $\left\{Q^j\right\}_{j\in{1,\dots,n}}$ such that
\begin{equation*}
\Lambda_{ij} :=  \innerprod{\mbf u_i}{\nabla^2\mbf{u}_j},\,\,\,
W_{ij} := -\innerprod{\mbf u_i}{S(\mbf u_j,\bar{\mbf u})},\,\,\,
Q^j_{ik} := -\innerprod{\mbf u_i}{\mbf u_j \cdot \nabla \mbf u_k},
\end{equation*}

with $L \eqdef \frac{1}{\Reynolds}\Lambda + W$, and defining a linear matrix-valued operator $N : \Re^{n} \to \Re^{n\times n}$ as 
\[
N(a) := a_jQ^j,
\]
one arrives at a compact representation of the ODE \eqref{eqn:nominalODE}
\begin{equation}\label{eqn:sysNominal}
\dot a = f(a) := La + N(a)a.
\end{equation}

Two useful general observations about this system are that the matrix $\Lambda$ is symmetric and negative-definite (in particular, it is diagonal if the basis functions $\mbf u_i$ are chosen as eigenfunctions of \eqref{eqn:EigenEnergy} with $e = 0$), and that $a\tpose N(a)a = 0$ for all $a$.  The latter assertion is a restatement of the energy conservation relation~\eqref{eqn:EnergyConservation} in finite dimensions.  

\section{Stability of Finite Dimensional Models using SOS}\label{sec:FiniteDim}

For simplicity of exposition, we will assume in this section that the steady solution $\bar{\mbf u}$ is spanned by the basis functions $\mbf u_i$, i.e.\ that there exist some real constants $c_i$ such that 
$\bar{\mbf u} = \mbf u_i c_i$.  In this case, one can rewrite the dynamics of the finite dimensional system~\eqref{eqn:sysNominal} in the equivalent form 
\begin{equation}\label{eqn:sysNominal:alt}
\dot a = \frac{1}{\Reynolds} \Lambda a + N(a + c)(a+c) - N(c)c,
\end{equation}
where we have used the identity $Wa = N(a)c + N(c)a$.  Note that $a = 0$ is an equilibrium solution to \eqref{eqn:sysNominal:alt}, and we wish to find the largest value of $\Reynolds$ for which this system is globally asymptotically stable.    

To this end, we first recall that an ODE system $\dot a = f(a)$ is stable if one can find a continuously differentiable Lyapunov function $V$ satisfying each of the following conditions~\cite[Thm 4.1]{Khal02}:
\begin{align}
V(0) &= 0 \tag{L1} \label{LyapConditions:1}\\
V(a) &>0 \quad \forall a \neq 0 \tag{L2} \label{LyapConditions:2}\\
\nablaa V(a) \cdot f(a) &< 0 \quad \forall a \neq 0. \tag{L3} \label{LyapConditions:3}
\end{align}
There is unfortunately no known method to construct such a function for an arbitrary system of nonlinear ODEs.  However, in the case of a system described exclusively by polynomial functions such as \eqref{eqn:sysNominal}, the situation is more hopeful.

First define the energy-like functions $E_ \theta : \Re^n \to \Re$ as 
\begin{align*}
E_\theta(a) &\eqdef \frac{1}{2} \norm{a + \theta c}^2.
\end{align*}
Of special interest will be the \emph{perturbation energy} function $E_0$ and the \emph{total energy} function~$E_1$.  
In particular, a useful observation is that  
\[
\nablaa E_0(a) \cdot N(a)a = a\tpose N(a)a = 0,
\]
 i.e.\ the nonlinear part of the dynamics of the system \eqref{eqn:sysNominal} is invariant with respect to the perturbation energy.  
Selecting as a candidate Lyapunov function $V = E_0$, 
stability of the system \eqref{eqn:sysNominal} is therefore assured for all $\Reynolds$ such that 
\begin{equation}\label{energyLyapCondition}
\nablaa V(a) \cdot f = a\tpose L a = a\tpose (\frac{1}{R}\Lambda + W) a < 0\quad \forall x \neq 0.
\end{equation}
Calculation of the maximum value of $\Reynolds$ for which~\eqref{energyLyapCondition} holds is then straightforward, since one needs only to find the largest $\Reynolds$ such that the matrix 
$2\Lambda + \Reynolds \cdot (W + W\tpose)$ remains negative-definite.  Of course, this mirrors exactly the situation in the infinite dimensional case.  

We next consider whether it is possible to establish stability of the system \eqref{eqn:sysNominal} using some alternative polynomial Lyapunov function.  In order to restrict the overall size of our search space, we first consider the essential features of such a function.

\subsection{System behavior for extreme values of $\norm {a}$}\label{sec:extremeNorm_a}

Consider first the linear part of the system~\eqref{eqn:sysNominal} in isolation, i.e.\ 
\begin{equation}\label{eqn:linearPart}
\dot a = L a.
\end{equation}
If $L$ has any positive eigenvalues then the system \eqref{eqn:linearPart} is unstable, implying immediately that the nonlinear system \eqref{eqn:sysNominal} is also unstable.  If the system \eqref{eqn:linearPart} is asymptotically stable, then
 there exists some $P \succ 0$ such that $V(a) = a\tpose P a$ is positive-definite and 
\begin{equation}\label{linearLyap}
\nablaa V(a) \cdot f(a) = a\tpose \left(L\tpose P + PL\right)a < 0 \quad \forall a \neq 0,
\end{equation}
see \cite[Thm\ 4.6]{Khal02}.
Such a function also ensures stability of the nonlinear system \eqref{eqn:sysNominal} for some region around the origin, since the linear component of \eqref{eqn:sysNominal} dominates when $\norm{a}\ll 1$. 

Considering the nonlinear term of~\eqref{eqn:sysNominal} in isolation, one typically expects that for any $V(a) = a\tpose Pa$ with $P \succ 0$,
\[
\set{a}{\nablaa V(a)  \cdot N(a)a > 0} \neq \emptyset,
\]
unless $P \propto I$.  Since the nonlinear component of~\eqref{eqn:sysNominal} is the dominant term when $\norm{a}\gg 1$, we should generally not expect to find a second-order positive-definite polynomial Lyapunov function $V$ other than the perturbation energy function $V = E_0$, or some monotone function thereof.  On the other hand, using the system representation \eqref{eqn:sysNominal:alt} it follows that
\begin{align}
\nablaa E_1(a) \cdot f(a)
&= (a + c) \cdot \left(\frac{1}{\Reynolds} \Lambda a + \left[N(a + c)\right](a + c)  -  N(c)c
				 \right)\notag\\[1ex]
&= \frac{1}{\Reynolds} \cdot \left(a\tpose \Lambda a + a\tpose \Lambda c\right) - a\tpose N(c)c.
                   \label{negdef_E1}
\end{align}

Consequently, $\nablaa E_1 \cdot f < 0$ for all $\norm{a}$ sufficiently large with respect to a fixed Reynolds number\footnote{
This effect is not exclusive to the total energy function $E_1$.  If one defines the energy-like function $E_d \eqdef \frac{1}{2}(a+d)\tpose (a+d)$, then the quadratic part of $
\nablaa E_d \cdot f$ is negative-definite whenever the Reynolds number $\Reynolds$ and vector $d$ are contained in the set
\[
\set{(\Reynolds, d)}{\Lambda + \frac{\Reynolds}{2} \left[ \Bigl(W + W\tpose\Bigr)  - \Bigl(\tilde W(d) + \tilde W\tpose(d) \Bigr)\right] \prec 0},
\]
where $\tilde W(d)$ is linear in $d$ and defined such that $ \tilde W(d)a = N(d)a + N(a)d$.  The above set is convex in $\Reynolds$ for fixed $d$ and vice-versa.  In the case that $\bar {\mbf u} = c_i{\mbf u_i}$, in follows that $W = \tilde W(c)$ and one 
can make the particularly convenient choice $d = c$, so that the above set is unbounded in $\Reynolds$.  
}, 
though the choice $V = E_1$ would not satisfy condition \eqref{LyapConditions:1}.
A reasonable approach therefore is to search for a candidate Lyapunov function in the form $V = A + B$, where the components $A : \Re^n \to \Re$ and $B : \Re^n \to \Re$ have the following properties:
\begin{subequations} \label{VeqAB}
\begin{alignat}{3}
[A + B](0) &= 0 & \label{VeqAB:}\\
\nablaa\left[ A(a) + B(a)\right] &\approx a\tpose P &\quad \forall \norm{a} &\ll 1 \label{VeqAB:1} \\
\nablaa B(a) &\approx \sigma(a) \cdot \nablaa E_1(a) &\quad \forall \norm{a} &\gg 1 \label{VeqAB:2} \\
\norm{\nablaa A(a)} &\ll \norm{\nablaa B(a)} &\quad\forall \norm{a} &\gg 1,
\label{VeqAB:3}
\end{alignat}
\end{subequations}
where $P \succ 0$ and $\sigma : \Re^n \to \Re$ is a nondecreasing positive function.
The condition \eqref{VeqAB:1} ensures that $\nablaa V \cdot f$ satisfies approximately the linear Lyapunov condition \eqref{linearLyap} in a localized region about the origin.  The conditions \eqref{VeqAB:2}--\eqref{VeqAB:3} ensure that $\nablaa V \cdot f < 0$ for all states sufficiently far from the origin, in accordance with \eqref{negdef_E1}.

In order to exploit SOS techniques, we restrict our attention to cases where both $A$ and $B$ are polynomial functions and $\deg A < \deg{B}$.  A useful observation is that any choice of $B$ in the form
\begin{equation}\label{EnergyProducts}
B(a) = \prod_{i=1}^k E_{\theta_i}(a) \,\,\, \text{with} \,\,\,\,\,
\frac{1}{k}\sum_{i=1}^k \theta_i = 1\,\,\, \text{and} \,\,\,\,\, \theta_1 = 0
\end{equation}
with $k \in \mbb{N}$
satisfies the condition \eqref{VeqAB:2}.  In searching for a Lyapunov function in the form $V = A+B$, we will view the function $A$ as a term to be optimized, and therefore refer to it as the \emph{variable term}.   We will restrict the function $B$ to be some combination of energy-like functions in the form~\eqref{EnergyProducts}, and hence refer to it as the \emph{energy term}.

\subsection{Lyapunov Function Generation Using Sum-of-Squares}

If we restrict our attention to polynomial functions $V$ with no constant term (so that $V(0) = 0$), then the Lyapunov conditions \eqref{LyapConditions:1}--\eqref{LyapConditions:3} can be rewritten as 
\begin{subequations}\label{Psatz:EmptySets}
\begin{align}
&\set{a}{
V(a) \le 0, \,\, \ell_1(a) \neq 0
} = \emptyset \\[2ex]
&\set{a}{-
\nablaa V(a)\cdot f(a) \le 0, \,\, \ell_2(a) \neq 0
} = \emptyset,
\end{align}
\end{subequations}
where positive-definite polynomial functions $\ell_i$ are used in place of the vector-valued condition $a \neq 0$. For simplicity, we can define the functions $\ell_i$ as
\[
\ell_i(a) \eqdef \sum_{j=1}^n \epsilon_{ij}a_j^2,
\]
and impose a strict positivity constraint on the values $\epsilon_{ij}$.  Straightforward application of the Positivstellensatz (see \cite{Parrilo00} and the references therein) shows that satisfaction of the conditions \eqref{Psatz:EmptySets} is assured if one can identify polynomial functions $(s_1,s_2)$ such that 
\begin{equation}\label{Psatz:Equalities}\tag{SOS}
\begin{aligned}
V(a) - \ell_1(a) &= s_1, \quad s_1 \in \Sigma_n 
\\
-\nablaa{V}(a) \cdot f(a)  - \ell_2(a) &= s_2, 
 \quad s_2 \in \Sigma_n, 
\end{aligned}
\end{equation}
where $\Sigma_n$ denotes the set of all sum-of-squares polynomials in $\Re^n$.  

The problem of determining whether \eqref{Psatz:Equalities} can be satisfied can be reformulated as a convex optimization problem in the form of a semidefinite program (SDP) using standard software tools \cite{yalmip,YalmipSOS09,PP05}.  If $\deg V = 2d$ (note that the degree of $V$ must be even for \eqref{LyapConditions:1} to be satisfied), then the general form of our problem is:
\begin{subequations}\label{SDP}
\begin{alignat}{3}
%
\hspace{-15ex}\smash{
	\text{
		\raisebox{-12ex}  
		{(SDP) $\left\{\rule{0pt}{13ex} 
		\right.$~~~~}
	}
}
%
%
&&\min_{V,H_1,H_2,\{\epsilon_{ij}\}} 
\,\,\,0 \,\,\,\quad\qquad\qquad\\
&&\text{subject to:}\qquad\qquad\qquad\quad\notag\\
&&V(x) - \ell_1(a) &= m_d(a)\tpose H_1 m_d(a) \label{SDP.H1}\\
&&-\frac{\partial V}{\partial a} \cdot f(a)  - \ell_2(a) &= m_d(a)\tpose H_2 m_d(a) \label{SDP.H2}\\
&&(H_1,H_2) &\succeq 0 \label{SDP_psd}\\
&&(\epsilon_{1,j},\epsilon_{2,j}) &\ge \bar\epsilon \quad \forall j \in \{1,\dots,n\}, \label{SDP_eps}
\end{alignat} 
\end{subequations}

where $m_d(a)$ is a vector of all monomials in $a$ with degree less than or equal to~$d$.  
The objective function in our optimization problem is zero since we are interested only in feasibility.  Note that any solution to the problem (SDP) will satisfy the original sum-of-squares condition \eqref{Psatz:Equalities}, since the semidefiniteness constraint \eqref{SDP_psd} ensures that \mbox{\eqref{SDP.H1}--\eqref{SDP.H2}} can be expressed as sums-of-squares following a suitable similarity transformation.  The lower bounding constant $\bar\epsilon$ for the terms $\epsilon_{ij}$ in \eqref{SDP_eps} must be strictly positive, though it is otherwise arbitrary.

\subsection{Determining Stable Values for $\Reynolds$}

One can estimate an upper bound on the value of $\Reynolds$ for which a solution to \eqref{Psatz:Equalities} can be found via a straightforward binary search strategy. However, in all but the trivial case $V = E_0$,  there is no reason to suppose a priori that if a solution to \eqref{Psatz:Equalities} can be found for some $\bar{\Reynolds}$, then a solution 
can be found for all 
$\Reynolds \in (0 \,\, \bar\Reynolds]$.  Provision of such an assurance is possible by augmenting \eqref{Psatz:Equalities} with additional constraints.  First note that the Lyapunov condition \eqref{LyapConditions:3} can be written as 
\begin{equation}\label{LyapConditions:3:exp}
\nablaa V(a) \cdot \biggl(W a + N(a)a\biggr) + 
\frac{1}{\Reynolds} \cdot \biggl(\nablaa V(a) \cdot \Lambda a \biggr)  < 0 \quad \forall a \neq 0.
\end{equation}
If \eqref{LyapConditions:3:exp} is satisfied for some $\bar\Reynolds$, then it is satisfied for all $\Reynolds \in [0,\bar\Reynolds]$ provided that the second term
$\nablaa V(a)\cdot \Lambda a \le 0$
for all $a$.  Satisfaction of this condition can be imposed as a sum-of-squares constraint,
\begin{equation}\label{SDP:lowReConstraint}
-\nablaa V(a) \cdot \Lambda a = s_3, \quad s_3 \in \Sigma_n,
\end{equation}
and included as an additional condition to \eqref{Psatz:Equalities} (or checked a posteriori).

Given a Lyapunov function $V$ for some value $\bar\Reynolds$, it is possible to compute directly the smallest and largest value $\Reynolds$ for which $V$ is a Lyapunov function, since \eqref{LyapConditions:3:exp} is affine in ${1}/{\Reynolds}$; e.g.\ one can compute an upper bound by solving the sum-of-squares problem
\begin{gather*}
\min_{\Reynolds \ge 0} \qquad{1}/{\Reynolds} \hspace{57ex}\\
\text{subject to:}\quad
-\nablaa V(a) \bigl(Wa + N(a)a\bigr) - 
\frac{1}{\Reynolds} \bigl(\nablaa V(a) \cdot \Lambda a \bigr)
- \ell_2(a) \in \Sigma_N 
\end{gather*} 
and taking the inverse of its minimum value.  

\subsection{Computational Complexity}

We next consider the computational effort required to solve the problem (SDP) for various degrees of Lyapunov candidate function $V$.  If we assume that $V$ is a polynomial function with arbitrary coefficients and $\deg V = 2d$, then the monomial vector $m(a)$ is composed of $D$ distinct monomial terms, where
$$D = \frac{(n+d)!}{n! d!}.$$ 
Standard results from semidefinite programming ensure that one can solve the problem (SDP) in $\mcf{O}(\sqrt{D})$ iterations
using a primal-dual interior point method%
\footnote{
More precisely, one can guarantee that a primal-dual interior point algorithm will reduce the duality gap of its solution iterate to a multiple $\epsilon$ of its original value within $\mcf{O}(\ln(1/\epsilon)\sqrt{D})$  iterations.  The reader is referred to \cite{Wolkowicz00,Todd01,Boyd04} and references therein for an overview of algorithms and complexity results for semidefinite programming. 
}, 
with each iteration requiring $\mcf{O}(D^3)$ operations.  In practice, it is generally the case that the number of iterations required to solve a semidefinite programming problem is roughly constant with respect to increasing problem size, so the computation time is determined almost entirely by the per-iteration computation cost.

The rapid increase in computational burden with increasing system dimension means that SOS methods are likely to be applicable for relatively low dimensional models only, even if one assumes that the considerable degree of problem-specific structure inherent in  \eqref{SDP} can be somehow exploited (e.g.\ using a structured approach such as \eqref{VeqAB}).  In particular, it is not advisable to attempt to estimate the maximum stable Reynolds number in the infinite dimensional Navier-Stokes system \eqref{eqn:NSpert} via solution of a succession of problems in the form \eqref{SDP} with increasing dimension.  We therefore require a more indirect approach, whereby the finite-dimensional techniques of this section can be extended to the infinite-dimensional system \eqref{eqn:NSpert} without excessive additional computation.  We propose such an approach in the remainder of the paper.

\section{Infinite Dimensional Flow Models}\label{sec:infDimProblems}

We now return to the general case where $\mbf u_s \neq 0$, which we will view as an uncertain forcing term in our ODE.  In this case substituting \eqref{eqn:ansatz} into \eqref{eqn:NSpert:S} and taking an inner product of both sides results in a model similar to the ODE \eqref{eqn:sysNominal}, but with additional perturbation terms in $\mathbf u_s$, i.e.\
\begin{equation}\label{eqn:sysDisturbed}
\dot a = f(a) + \Theta_{a}(\mbf u_s) + \Theta_b(\mbf u_s,a) + \Theta_c(\mbf u_s),
\end{equation}
where the additional perturbation terms are defined as
\begin{subequations}\label{eqn:disturbanceTerms}
\begin{align}
\Theta_{ai}(\mbf u_s) &\eqdef \hphantom{+}
\frac{1}{\Reynolds}\innerprod{\mbf u_i}{\nabla^2\mbf{u}_s}  
-\innerprod{\mbf u_i}{S(\mbf u_s,\bar{\mbf u})} 
\\
\Theta_{bi}(\mbf u_s,a) &\eqdef - \innerprod{\mbf u_i}{S(\mbf u_j,\mbf u_s)}a_j\\
\Theta_{ci}(\mbf u_s) &\eqdef -\innerprod{\mbf u_i}{\mbf u_s \cdot \nabla \mbf u_s},
\end{align}
\end{subequations}
and $f(a)$ is as defined in \eqref{eqn:sysNominal}.
In the above, a subscript $i$ indicates that the expression is the $i$th element of a vector quantity.
The perturbation term $\Theta_a$ represents a linear disturbance in $\mbf{u}_s$, $\Theta_b$ represents a bilinear disturbance in $(\mbf{u}_s,a)$, and $\Theta_c$ represents a quadratic disturbance in $\mbf{u}_s$.    

We would like to bound the influence of each of these perturbation terms on our ODE in terms of $\norm{\mbf u_s}$ and $\norm{a}$. In order to do so, we apply \eqref{eqn:tripleproduct} repeatedly to eliminate the appearance of terms $\nabla \mbf u_s$, so that \eqref{eqn:disturbanceTerms} can be rewritten as%
\footnote{The notation used can be clarified by the equivalent expression for the Cartesian components of the vector $\mbf h_i$: $h_i^m= 
\frac{1}{\Reynolds}{\nabla^2 u_i^m}  
-\frac{\partial \bar{u}^k}{\partial x_m}  u_i^k 
+\bar{\mbf u} \cdot \nabla u_i^m.$
}
\begin{subequations}\label{eqn:disturbanceTerms:2}
\begin{align}
\hspace{13ex}  
\Theta_{ai}(\mbf u_s) 
&=\innerprod{\mbf u_s}{\mbf h_i}, & 
\mbf h_i&\eqdef 
\frac{1}{\Reynolds}{\nabla^2\mbf{u}_i}  
-\nabla \bar{\mbf u} \cdot \mbf u_i 
+\bar{\mbf u} \cdot \nabla \mbf u_i, 
\hspace{6ex}   
\\
\Theta_{bi}(\mbf u_s,a) 
&=\innerprod{\mbf u_s}{\mbf h_{ij}}a_j, &
\mbf h_{ij}&\eqdef\mbf u_j\cdot \nabla\mbf u_i -
\nabla\mbf u_j\cdot \mbf u_i, \\
\Theta_{ci}(\mbf u_s) 
&= 
\innerprod{\mbf u_s}{\mbf u_s \cdot \nabla \mbf u_i}.&
\end{align}
\end{subequations}

We are of course left with an ODE in the form \eqref{eqn:sysDisturbed} which still features the perturbations~$\mbf u_s$.  We next bound the influence of this term by modeling only the evolution of its energy $q$, which we model as $q^2 = \norm{\mbf u_s}^2/2$.  In the process we add a single ODE to supplement \eqref{eqn:sysDisturbed}, representing the time evolution of the squared energy term $q^2$.  

Substituting \eqref{eqn:ansatz} into \eqref{eqn:NSpert:S} and taking an inner product of both sides with the total velocity field $\mbf u = \mbf u_ia_i + \mbf u_s$ provides the additional ODE in term of the perturbation energy~$q^2$, 
\begin{align}
\dot{\left(q^2\right)}&= 
a\tpose f(a) - a\tpose \dot a + \Gamma(\mbf{u}_s) + \chi(\mbf u_s,a),\label{eqn:energyDissipation}
\shortintertext{where }
\Gamma(\mbf{u}_s) 
& \eqdef
\frac{1}{\Reynolds}\innerprod{\mbf u_s}{\nabla^2\mbf u_s}
- \innerprod{\mbf u_s}{S(\mbf u_s,\bar{\mbf u})},\label{eqn:Gamma}\\
\chi(\mbf u_s,a) & \eqdef \innerprod{\mbf u_s}{\mbf g_j}a_j,\quad \mbf g_j \eqdef 
\left(\frac{1}{\Reynolds}\nabla^2- e\right) \mbf u_j.
\end{align}

Verification of the above relies on the aforementioned assumptions about the subspace $\mcf{S}$ and on application of the various identities described in Section \ref{sec:intro}.  In particular, these allow one to establish the relations
\begin{align*}
\innerprod{\mbf u}{\frac{\partial \mbf u}{\partial t}} = a\tpose \dot a + \dot{(q^2)} \text{~~and~~}
a\tpose f(a)
= \left[\frac{1}{\Reynolds}\innerprod{\mbf u_i}{\nabla^2 \mbf u_j} - \innerprod{\mbf u_i}{S(\mbf u_j,\bar{\mbf u})}\right]a_ia_j.
\end{align*}

Note that in~\eqref{eqn:energyDissipation}, the terms $a\tpose f(a)$ and $\Gamma(\mbf u_s)$ represent the self-contained dissipation or generation of energy depending on $a_i\mbf u_i$ and $\mbf u_s$, while the term $\chi(\mbf u_s, a)$ represents the generation or dissipation of energy containing cross terms between these velocity fields.

\subsection{Description as an Uncertain System}
\label{DescriptionasanUncertainSystem}
The complete system of interest can now be written as%
\begin{subequations}\label{eqn:completeSystem}
\begin{align}
\dot a &= f(a)
 + \Theta_a(\mbf u_s) + \Theta_b(\mbf u_s,a) + \Theta_c(\mbf u_s)
\label{eqn:completeSystem:a}\\
\dot{\left(q^2\right)}
&= 
a\tpose f(a) - a\tpose \dot a + \Gamma(\mbf{u}_s) + \chi(\mbf u_s,a).\label{eqn:completeSystem:b}
\end{align}
\end{subequations}

We are now free to treat $\mbf u_s$ as an uncertain term driving the ODE system \eqref{eqn:completeSystem}, whose time evolution is known to satisfy the subspace constraint $\mbf u_s \in \mcf {S}$ and the energy constraint~\mbox{$q^2=\norm{\mbf u_s}^2/2$}.  The worst-case effect of this uncertainty can then be bounded via appropriate norm bounds.  

The first of these bounds relates to the uncertain terms in \eqref{eqn:completeSystem:a}.  
There exist constants $c_i \ge 0$ and a polynomial function $p_1(a,q) \ge 0$ such that 
\begin{equation}\label{lem:normSquareBounds}
\norm{\Theta_{a}(\mbf u_s) + \Theta_b(\mbf u_s,a) + \Theta_c(\mbf u_s)}^2 
\le p_1(a,q)=
c_1 q^2 + c_2 q^2\norm{a}^2 + c_3 q^4 
\end{equation}
for any $a$ and $\mbf{u}_s.$  A rigorous proof of the existence of these constants is given in~\ref{app:NormBoundConstants}. Critically, estimation of the coefficients $c_i$ involves the solution only of linear problems for partial differential equations and optimization over finite-dimensional polynomials.

A second bound relates to the uncertain term $\Gamma(\mbf u_s)$ in \eqref{eqn:completeSystem:b}.  Comparing \eqref{eqn:Gamma} with \eqref{eqn:EnergyInequality} shows that $\Gamma(\mbf{u}_s)\le\kappa \norm{\mbf u_s}^2$ with $\kappa = \lambda_1$.   However, $\mbf u_s$ satisfies the additional constraints $\innerprod{\mbf u_s}{\mbf u_i}=0$ and therefore may admit a stronger bound. Note that the number of positive eigenvalues of \eqref{eqn:EigenEnergy} is always finite \cite{Batcho2001}. Hence, if $\mbf u_i$ are chosen as the first $n$ eigenfunctions of \eqref{eqn:EigenEnergy} and $n$ is large enough, then 
\begin{equation}
\label{assum:us:energyStable}
\Gamma(\mbf u_s) \le \kappa_s \norm{\mbf{u}_s}^2 = 2\kappa_s q^2
\end{equation}
for all $\mbf u_s \in \mcf{S}$, where $\kappa_s=\lambda_{n+1} < 0.$  
If $\mbf u_i$ are not eigenfunctions of \eqref{eqn:EigenEnergy}, then $\kappa_s$ is the largest eigenvalue of the following problem
\[
\begin{gathered}
\lambda \mbf u + \mu_k \mbf u_k  = (e - \frac{1}{\Reynolds}\Delta)\mbf u + 
\nabla p \\
\innerprod{\mbf u_k}{ \mbf u} = 0,\,\,
\nabla \cdot \mbf u = 0,\,\, \left.\mbf u\right|_{\partial\Omega}=0.
\end{gathered}
\]
In what follows we will assume that $\kappa_s <0$ in \eqref{assum:us:energyStable}. 

A final bound relates to the uncertain term $\chi(\mbf u_s)$ in \eqref{eqn:completeSystem:b}. 
If $\mbf u_i$ are eigenfunctions of \eqref{eqn:EigenEnergy} then $\chi=0,$ because in this case $\mbf g_i=-\lambda_i\mbf u_i+\nabla\phi_i$ with some scalar functions $\phi_i$ and because $\mbf u_s$ is orthogonal to both $\mbf u_i$ (by definition) and to gradients of any scalars (since $\nabla \cdot \mbf u_s=0$). In the general case 
there exists a constant $d$ and a polynomial function $p_2(a,q) \ge 0$ such that
\begin{equation}\label{lem:normChiBound}
\norm{\chi(\mbf u_s,a)}^2 \le p_2(a,q)=d q^2\norm{a}^2.
\end{equation}
The proof is very similar to the proof of \eqref{lem:normSquareBounds}.

\section{Stability of Infinite Dimensional Models using SOS}\label{sec:InfiniteDim}

Given the (uncertain) ODE system \eqref{eqn:completeSystem}, we can now search for a Lyapunov function verifying stability of the composite state vector $(a,q^2)$.  We therefore would like to construct a Lyapunov function $V : \Re^n \times \Re \to \Re$ such that
\begin{equation}\label{eqn:lyapCondition}
\frac{\partial V}{\partial a} \dot a + \frac{\partial V}{\partial(q^2)} \cdot  \dot{\left(q^2\right)} < 0, \quad \forall (a,\mbf u_s) \ne 0, \,\,\mbf u_s \in S.
\end{equation}

We can expand the left hand side of this condition and collect terms to get the equivalent Lyapunov condition
\begin{equation} \label{eqn:lyapCond}
\frac{\partial V}{\partial a} f + 
\frac{\partial V}{\partial(q^2)} \Gamma + 
\left(\frac{\partial V}{\partial a} - \frac{\partial V}{\partial(q^2)}a\tpose\right)
\biggl(\Theta_
{a} + \Theta_b + \Theta_c\biggr) + \frac{\partial V}{\partial(q^2)} \chi < 0,
\end{equation}
where we have omitted the arguments for $(f,\Gamma,\Theta_a,\Theta_b,\Theta_c,\chi)$ for brevity.  
For simplicity, we will assume that the function $V$ is chosen in such a way that 
\begin{equation}\label{assum:dVdq2}
{\partial V}/{\partial(q^2)}\ge0
\end{equation}
Consequently, 
\[
\frac{\partial V}{\partial(q^2)} \cdot  \Gamma(\mbf{u}_s) < 0, \quad  \forall \mbf u_s \in \mcf{S} \backslash \{0\}
\]
provided that \eqref{assum:us:energyStable} and~\eqref{assum:dVdq2} hold.

\subsection{Comparison to the Energy Method}\label{sec:energyCompare}

Note that if one chooses a candidate Lyapunov function by making the most obvious generalization of the type of function suggested in Section \ref{sec:FiniteDim}, i.e.\ if one chooses 
$
V(a,q^2) = V_a(a) + \prod_{i=1}^{k}(E_{\theta_i}(a)  + q^2),
$ 
where $V_a(\cdot)$ is some polynomial function, then \eqref{assum:dVdq2} is satisfied. 

The term 
\[
\left(\frac{\partial V}{\partial a} - \frac{\partial V}{\partial(q^2)}a^\top\right)
\]
in the Lyapunov condition~\eqref{eqn:lyapCond} can be viewed as a misalignment between the (scaled) gradient of the energy function $E_0(a)$ and the gradient term ${\partial V}/{\partial a}$. If one chooses as a candidate Lyapunov function
\[
V(a,q^2) =E_0(a) + q^2,
\]
then the above misalignment term is zero.  If, additionally, one chooses $\mbf u_i$ such that $\chi=0,$ which was shown above always to be possible,  the situation reduces to the usual global stability condition using energy functions.  
  Consequently, if energy can be used as a Lyapunov function for the system \eqref{eqn:NSpert} for some Reynolds number $\Reynolds$, then the choice $V(a,q^2) = E_0(a) + q^2$ will satisfy the conditions~\eqref{eqn:lyapCond}. 

When the system remains globally stable for Reynolds numbers beyond this energy stability limit, one should first ask whether there exists \emph{any} polynomial in $(a,q^2)$ that will serve as a Lyapunov function.  We give a constructive proof of the existence of such a function in \ref{app:existence}.

It remains to demonstrate that the Lyapunov function $V$ satisfying~\eqref{eqn:lyapCondition} can be constructed in a systematic way using the SOS approach.

\subsection{Conversion to a Sum-of-Squares problem}\label{ConversiontoaSumofSquaresproblem}

After applying \eqref{assum:us:energyStable}, the inequality in the Lyapunov condition \eqref{eqn:lyapCond} can be written in vectorized form as
\begin{align}\label{eqn:lyapCond_approx}
\left(\frac{\partial V}{\partial a} f + 
\frac{\partial V}{\partial(q^2)} \cdot 2\kappa_s q^2 \right)
< -
\begin{bmatrix}
\left( 
\frac{\partial V}{\partial a} - \frac{\partial V}{\partial(q^2)}a\tpose 
\right), & \frac{\partial V}{\partial(q^2)}
\end{bmatrix}
\begin{bmatrix}\Theta_{a} + \Theta_b + \Theta_c\\
\chi
\end{bmatrix}.
\end{align}
We next apply the Schwarz inequality, \eqref{lem:normSquareBounds}, \eqref{lem:normChiBound}, and \eqref{assum:dVdq2} to arrive at a sufficient condition for satisfaction of the inequality \eqref{eqn:lyapCondition}:
\begin{equation} \label{eqn:lyapCond_Schur}
\!\left(\frac{\partial V}{\partial a} f(a) \!+\! 
\frac{\partial V}{\partial(q^2)} \cdot 2\kappa_s q^2 \right)
\!<\! \!-\! 
\left|
\begin{bmatrix}
\left(
\frac{\partial V}{\partial a} \!-\! \frac{\partial V}{\partial(q^2)}a\tpose 
\right), \!& \frac{\partial V}{\partial(q^2)}
\end{bmatrix}
\right|
\left[p_1(a,q) + p_2(a,q)\right]^\half,
\,\, \forall(a,q) \neq 0,
\end{equation}

where $|\cdot |$ represents the standard  Euclidian norm in $\Re^n$.  The above can be rewritten  more compactly as 
\begin{equation}\label{eqn:lyapCond_Schur_compact}
g(a,q) < -\norm{h(a,q)} \cdot  p^{\half}(a,q),\quad \forall(a,q) \neq 0,
\end{equation}
where
\begin{align*}
g(a,q) &\eqdef \frac{\partial V}{\partial a} f(a) + 
\frac{\partial V}{\partial(q^2)} \cdot 2\kappa_s q^2 \\
h(a,q) &\eqdef \begin{bmatrix}
\left(
\frac{\partial V}{\partial a} - \frac{\partial V}{\partial(q^2)}a\tpose 
\right), & \frac{\partial V}{\partial(q^2)}
\end{bmatrix}\tpose\\
p(a,q) &\eqdef p_1(a,q) + p_2(a,q).
\end{align*}

We next apply the following matrix property, based on the Schur complement \cite[Sec.~A.5.5]{Boyd04}.  For any vector $u$ and scalar $t$, 
\[
\norm{u} < t \Longleftrightarrow
\begin{bmatrix}
t & u\tpose \\ u & tI
\end{bmatrix} \succ 0.
\]
The condition \eqref{eqn:lyapCond_approx} is therefore equivalent to
\[
\begin{bmatrix}
g(a,q) & h\tpose(a,q)\cdot  p^\half(a,q) 
\\
h(a,q)\cdot  p^\half(a,q) 
&
g(a,q) I
\end{bmatrix} \prec 0, \quad \forall (a,q)\neq 0,
\]

Since $ p$ is a nonzero function for all nonzero $(a,q)$, the above can be pre- and post-multiplied by $\left(\begin{smallmatrix} p^\half & 0 \\ 0& I\end{smallmatrix}\right)$ to get the equivalent matrix inequality condition
\begin{equation}\label{eqn:lyapLMI}
 H(a,q) \eqdef
\begin{bmatrix}
 g(a,q)\cdot p(a,q) & h\tpose(a,q)\cdot  p(a,q) 
\\
h(a,q)\cdot  p(a,q) 
&
g(a,q) I
\end{bmatrix} \prec 0, \quad \forall (a,q)\neq 0.
\end{equation}

The most important thing to note about \eqref{eqn:lyapLMI} is that it is \emph{linear} in the coefficients of the Lyapunov function $V$.  This linear matrix inequality~(LMI) can be converted to an equivalent {scalar} polynomial inequality via introduction of an additional variable $z \in \Re^{n+1}$.  It is straightforward to verify that~\eqref{eqn:lyapLMI} is equivalent to
 \begin{equation}\label{eqn:lyapLMI_scalarized2}
z\tpose  H(a,q)z
< 0, \quad \quad\forall z \neq 0, \, \forall(a,q) \neq 0.
\end{equation}
%


%
In other words, if the function $V$ is chosen such that ${\partial V}/{\partial (q^2)}$ is nonnegative and \\[-0.2em]%
\begin{subequations}\label{eqn:Psatz_b}
\begin{gather}
\set{(a,q,z)}{
\begin{gathered}
z\tpose  H(a,q) z \ge 0\\
\norm{a}^2 + q^2 \neq 0,\,\, \norm{z}^2 \neq 0
\end{gathered} 
}
= \emptyset \label{eqn:Psatz1_b},
\intertext{and}
\set{(a,q)}{
\begin{gathered}
V(a,q^2) \le 0\\
\norm{a}^2 + q^2 \neq 0
\end{gathered} 
}
= \emptyset \label{eqn:Psatz2_b}
\end{gather}
\end{subequations}
then the Lyapunov condition \eqref{eqn:lyapCond} is satisfied.  
This is a standard form convenient for applying the Positivstellensatz theorem~\cite{Parrilo2003}, from which it follows that \eqref{eqn:Psatz_b} are satisfied if and only if there exist non-negative integer values $M_1, M_2, M_3, M_4$ and $M_5$ and sum-of-squares of polynomials $s_i(a,q,z)$ and $\sigma_j(a,q)$ such that
$$
s_0=-\sum_{i=1}^{M_1} s_i [z\tpose  H(a,q) z)]^i-(\norm{a}^2+q^2)^{2M_2}\norm{z}^{4M_3}
$$
and 
$$
\sigma_0=-\sum_{j=1}^{M_4}\sigma_j [V(a,q^2)]^j-(\norm{a}^2+q^2)^{2M_5}.
$$ 
One can now choose (from empirical considerations or by trial and error) the integers $M_j$ and the polynomials $s_i$ and $\sigma_j$ for $(i,j)>0,$ thus obtaining expressions for $s_0$ and $\sigma_0$ via the coefficients of the polynomial $V.$ Then determining the coefficients of $V$ such that $s_0$ and $\sigma_0$ are sum of squares can be attempted using the existing packages SOSTOOLS \cite{PP05} or YALMIP \cite{yalmip,YalmipSOS09}. For the purposes of computational efficiency, the selection of $V$ should be subject to structural constraints similar to the case of an ODE system obtained by simple truncation, as in Section \ref{sec:FiniteDim}.  

Finally, we note that if one is able to choose the bases $\mbf u_i$ such that $\chi = 0$, then the preceding problem can be simplified somewhat.  In this case, one is free to define 
\begin{align*}
h(a,q) &\eqdef
\left(
\frac{\partial V}{\partial a} - \frac{\partial V}{\partial(q^2)}a\tpose 
\right)\tpose, \quad
p(a,q) \eqdef p_1(a,q),
\end{align*}
in \eqref{eqn:lyapCond_approx}--\eqref{eqn:lyapCond_Schur_compact}.  It is then easy to show that the choice of Lyapunov function $V = E_0(a) + q^2$ will satisfy not only the stability condition~\eqref{eqn:lyapCond}, but also the robust LMI condition \eqref{eqn:lyapLMI} whenever the total perturbation energy is a Lyapunov functional for the Navier-Stokes system \eqref{eqn:NSpert}.  This ensures that the proposed method will always yields results at least as good as classical energy perturbation methods.  

\section{A Finite Dimensional Example}\label{sec:example}
\begin{figure}[h]\centering
{\includegraphics[scale=0.6,draft=\isDraft]{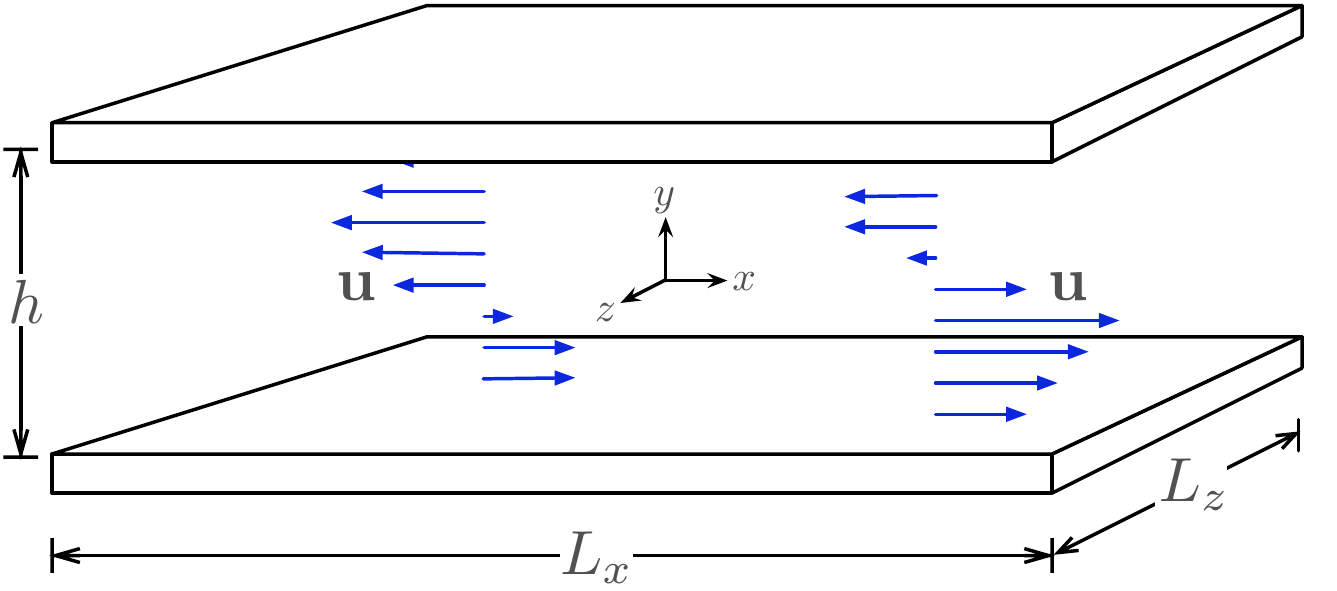}}
\caption{Periodic flow between parallel plates.}
\label{fig:CouetteFlow}
\end{figure}

In this section we present numerical results for a model of Couette flow\footnote{The results of this section were first reported in \cite{GC10}.} using the finite-dimensional stability analysis results of Section \ref{sec:FiniteDim}.  Couette flow refers to the shear flow of a fluid between two infinite parallel plates as shown \mbox{in Figure~\ref{fig:CouetteFlow}}.  

We employ the finite-dimensional ninth-order model developed in~\cite{MFE04} for this flow, and make assumptions identical to those in \cite{MFE04,MFE05} for the purposes of comparison.  The volume force is assumed to be
\[
\mbf{f} = \Bigl(\frac{\sqrt{2}\pi^2}{4\Reynolds}\sin(\pi y/2);0;0\Bigr),
\]
and the flow is assumed to be periodic in the spanwise and streamwise directions, with wavelengths $L_z$ and $L_x$ respectively.  We fix $L_x = 4\pi$ and $L_z = 2\pi$ throughout, with plate separation $h = 2$.

The flow is assumed to have free-slip boundary conditions
\[
u_y |_{y = \pm 1} = 0, \quad 
\left.\frac{\partial u_x}{\partial y} \right|_{y = \pm 1} = 
\left.\frac{\partial u_z}{\partial y} \right|_{y = \pm 1} = 0.
\]
A set of nine basis functions $\mathbf{u}_i$ were selected in~\cite{MFE04} based on physical insights and observations  arising from numerical simulation and experiment.  For reference these basis functions and their expansion into a nonlinear ODE in the form \eqref{eqn:sysNominal:alt} (equivalently \eqref{eqn:sysNominal}) are included in \ref{app:example}.  In this example, $\mathbf{u}_1$ is a laminar solution to the Navier-Stokes equation~\eqref{eqn:NS}. 

We compute an upper bound on the value of $\Reynolds$ for which the system~\eqref{eqn:sysNominal} is guaranteed to be stable using a variety of Lyapunov candidates, each of which satisfies the structural conditions~\eqref{VeqAB}--\eqref{EnergyProducts}.  

To compute the upper bound, we use a bisection method to find the largest Reynolds number for which the sum-of-squares optimization problem \eqref{SDP} could be solved with \mbox{$\bar\epsilon = 10^{-5}$} in \eqref{SDP_eps}.   All of the results obtained were computed on a 2.33 GHz Intel Xeon processor with 3.6 GB RAM, using the YALMIP interface to the SDP solver SeDuMi~\cite{yalmip,sedumi}.

Overall results are summarized in Table \ref{tab:LyapResults}.  
In each case we provide the form of Lyapunov function used and the maximum value of $\Reynolds$ for which a Lyapunov function in this form could be identified.   We also report the total solver time required (which includes the time spent both in the SDP solver SeDuMi and in preprocessing tasks by YALMIP), the total number of monomial terms in the vector $m_d(a)$ that appears in the equality constraints of \eqref{SDP.H1}--\eqref{SDP.H2}, the number of decision variables in $V$ that take nonzero values in the solution to $\eqref{SDP}$, and the total number of nonzero elements required in the solution for the matrices $H_1$ and $H_2$.  

 The largest Reynolds number for which a Lyapunov function could be identified was $\Reynolds = 54.1$, which compares very favorably to the value $\Reynolds = 7.5$ for which a purely energy-based method succeeds. 
 Previous results from numerical work in 
\cite{MFE05}, using the same model, have suggested a value of no more than $\Reynolds \simeq 80$ before the system becomes unstable.  This suggests that the method we propose is not unduly conservative.    Note that all of the computed stability bounds on $\Reynolds$ presented in Table~\ref{tab:LyapResults} are unchanged if one includes the additional constraint \eqref{SDP:lowReConstraint}.

\newcommand{\ivs}{\rule{0pt}{3ex}}
{
\begin{table}[ht]
\footnotesize
\begin{tabular}{|c|l|c|c|c|c|c|}\hline
\rule{0pt}{4ex}
{Case}& 
Lyapunov Function  & 
{$\Reynolds_{\max}$} & 
{\parbox[c]{1.55cm}{\center solver \\time (sec)}} & 
{\parbox[c]{1.55cm}{number of\\ monomials}}  & 
{\parbox[c]{1.55cm}{$\nnz(P)$ or \\\hphantom{~~~}$\dim(p)$}} & 
{\parbox[c]{1.55cm}{$\nnz(H_1)+$\\\hphantom{~~~}$\nnz(H_2)$}}
\\[2ex]\hline
%
\ivs 1&$V = E_0$ & 7.5  & 0.1 & -- & -- & -- \\\hline
%
\ivs 2&$V = \frac{1}{2}a\tpose P a + E_0E_2$  & 23.9   & 2.7 
& 54 & 21 & 1276\\ \hline
%
\ivs 3&$V = \frac{1}{2}a\tpose P a + E_0E_1E_2$  & 28.5   & 20.8 & 219 & 21 & 19410 \\ \hline
%
\ivs 4&$V = m_2(a)\tpose P m_2(a) + E_0E_1E_2$  & 54.1   & 43.3 & 219 & 776 & 24042 \\ \hline
%
\ivs 5&$V = p\tpose m_4(a) + E_0E_1E_2$  & 54.1   & 41.4  & 219 & 190 & 24042  \\ \hline
\end{tabular}
\caption{
Computed lower bounds on maximum stable Reynolds number for various Lyapunov functions \vspace{-4ex}}
\label{tab:LyapResults}
\end{table}
}
\vspace{-1ex}
\subsection{Perturbation Energy as a Lyapunov Function [Case 1]}\label{sec:Numerics:energyV}

We first consider the simplest case where one takes $V = E_0$.  Recalling \eqref{energyLyapCondition}, an upper bound on the value of $\Reynolds$ for which the system is guaranteed stable is readily found via solution of the following semidefinite programming (SDP) problem:
\begin{align*}
\max & ~~~\Reynolds \\
\textrm{s.t.~}&{~~~ 2 \Lambda + \Reynolds \!\cdot\!  (W + W\tpose)}  \prec 0.
\end{align*}
This method is analogous to the use of the conventional energy-based approach described in \cite{Serrin59,Joseph76}.

\subsection{Lyapunov Functions with Second-Order Variable Terms [Cases 2 \& 3]}
We next consider candidate Lyapunov functions with the variable term $A$ restricted to the quadratic form \mbox{$A(a) = a\tpose P a$}, where $P \in \mbb S^n$ is taken as a decision variable to be optimized, and $B$ restricted to a weighted sum of energy terms of higher order.  We restrict $P$ to those matrices whose sparsity patterns match that of solutions to the \mbox{Lyapunov~equation}
\begin{equation}\label{ex:Psparsity}
L \tpose P + PL = -I,
\end{equation}

and note that the sparsity pattern of such solutions are invariant with respect to $\Reynolds$.
%

This choice does not result in any apparent increase in conservatism (i.e.\ it does not affect the range of values for which \eqref{Psatz:Equalities} can be solved), while reducing considerably the overall computation time required.  Increasing the degree of the energy term between cases 2 and 3 shows a slight improvement in the maximum value of $\Reynolds$ for which stability can be assured.

\subsection{Lyapunov Functions with Fourth-Order Variable Terms [Cases 4 \& 5]}

Finally, we consider candidate Lyapunov functions with variable terms of fourth order.  In these cases, direct solution of \eqref{Psatz:Equalities} requires substantially increased computational effort relative to cases with second-order variable terms.  For case 4 we take $A(a) = m_2(a)\tpose P m_2(a)$ with the matrix~$P$~treated as a symmetric decision variable.  For case~5 we take~\mbox{$A(a) = p\tpose m_4(a)$} with vector $p$ treated as a decision variable.  In both cases, $A$ is restricted to contain only terms of degree at least~2.

In both cases, a solution to \eqref{Psatz:Equalities} was found once initially for an arbitrarily chosen (small) value of $\Reynolds$ in order to identify a likely sparsity pattern for $P$ and $p$ respectively. Subsequent computations enforced sparsity of $P$ and $p$ by setting to zero those elements taking relatively small values (i.e.~$\le 10^{-7}$) in the first trial.  In both cases, this procedure resulted in a substantial reduction in the degrees of freedom afforded to the solver, with a consequently large reduction in overall computation time.  

\begin{figure}[!htpb]
\newcommand{\polarFigWidth}{0.5\columnwidth}
\newcommand{\polarFigHeight}{0.48\columnwidth}
\newcommand{\polarFigBoxSize}{0.48\columnwidth}
~\\[1.5ex]
\centering
\makebox[\polarFigBoxSize]{
    \subfigure[Cut through $(a_1,a_3)$.]%
    {
      \label{phasePlot:13}
      \!\!\includegraphics[width=\polarFigWidth,height=\polarFigHeight,draft=\isDraft]
      				  {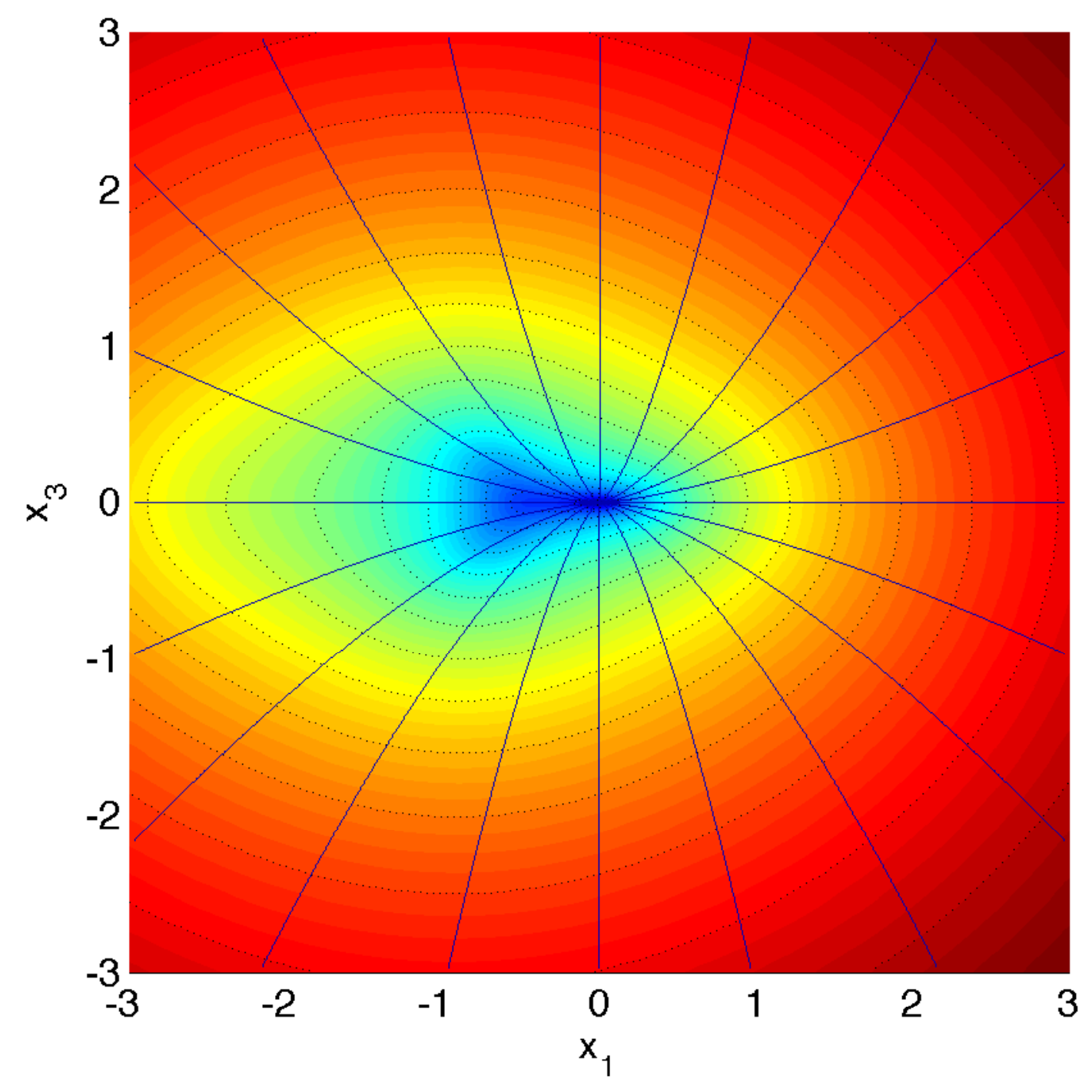}%
    }%
}%
\hfill
\makebox[\polarFigBoxSize]{
    \subfigure[Cut through $(a_6,a_7)$.]%
    {
      \label{phasePlot:67}
      \includegraphics[width=\polarFigWidth,height=\polarFigHeight,draft=\isDraft]
      				  {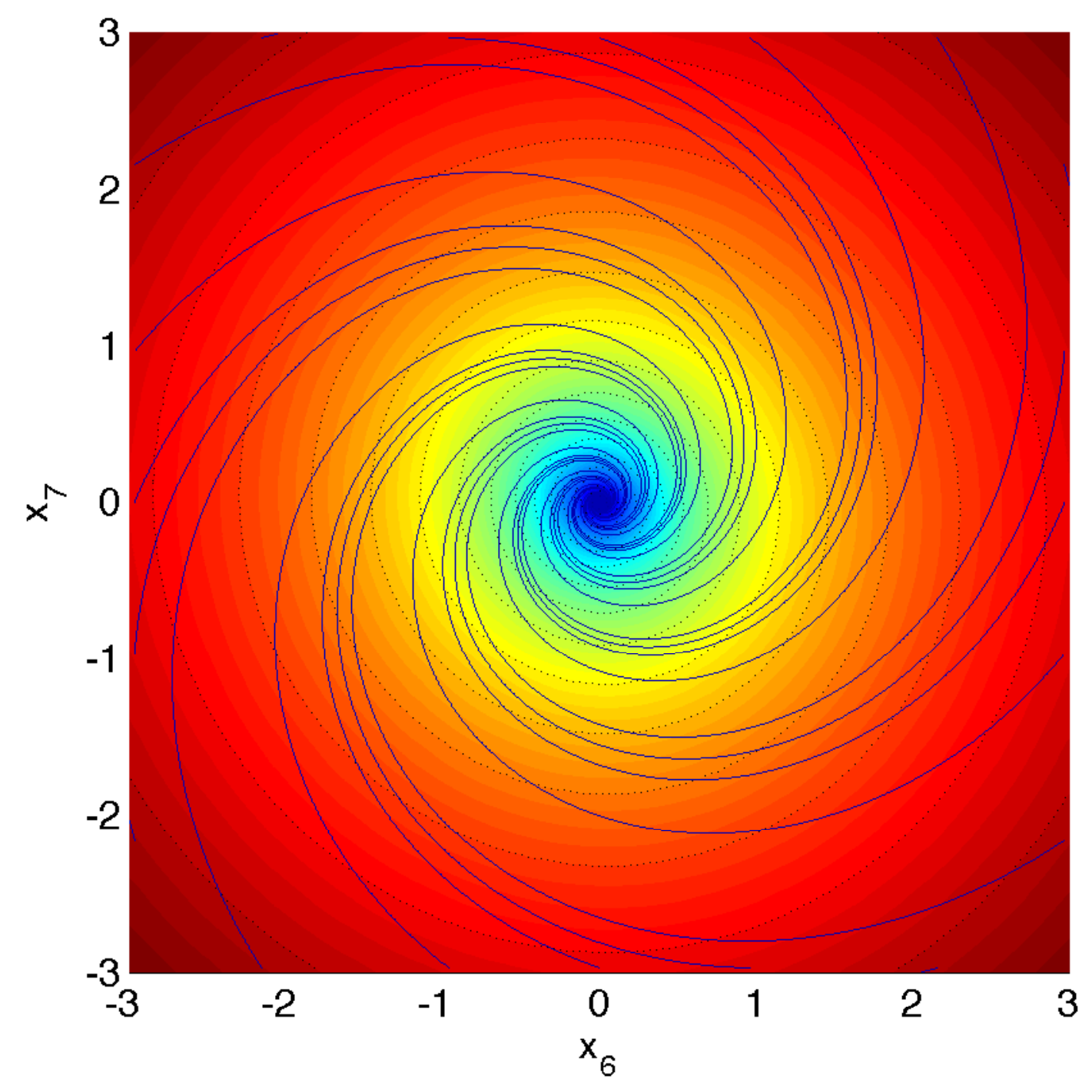}%
    }%
}~~\\
\makebox[\polarFigBoxSize]{
    \subfigure[Cut through $(a_1,a_9)$.]%
    {
      \label{phasePlot:19}
      \!\!\includegraphics[width=\polarFigWidth,height=\polarFigHeight,draft=\isDraft]
      				  {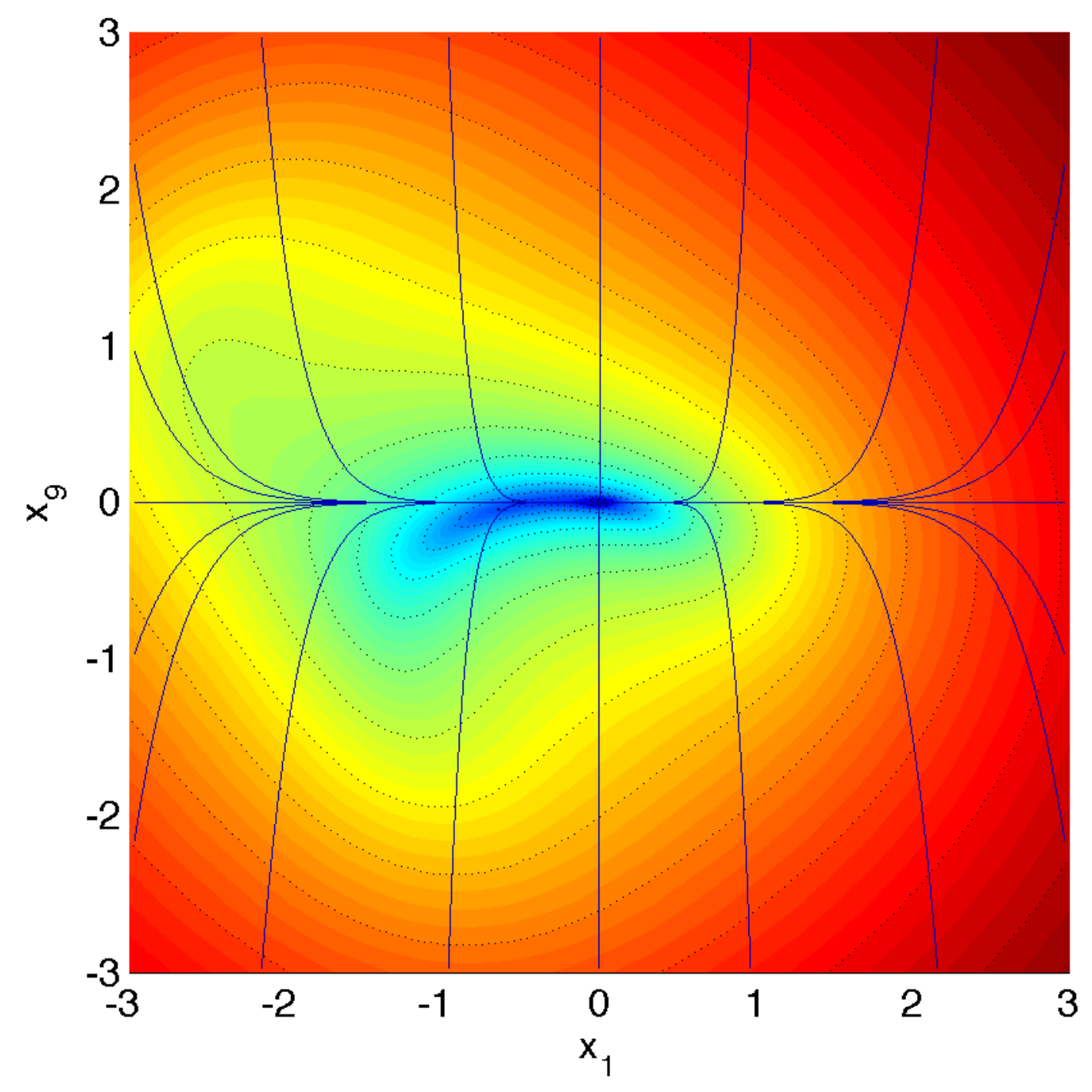}%
    }%
}%
\hfill
\makebox[\polarFigBoxSize]{
    \subfigure[Cut through $(a_2,a_9)$]%
    {
      \label{phasePlot:29}
      \includegraphics[width=\polarFigWidth,height=\polarFigHeight,draft=\isDraft]
      				  {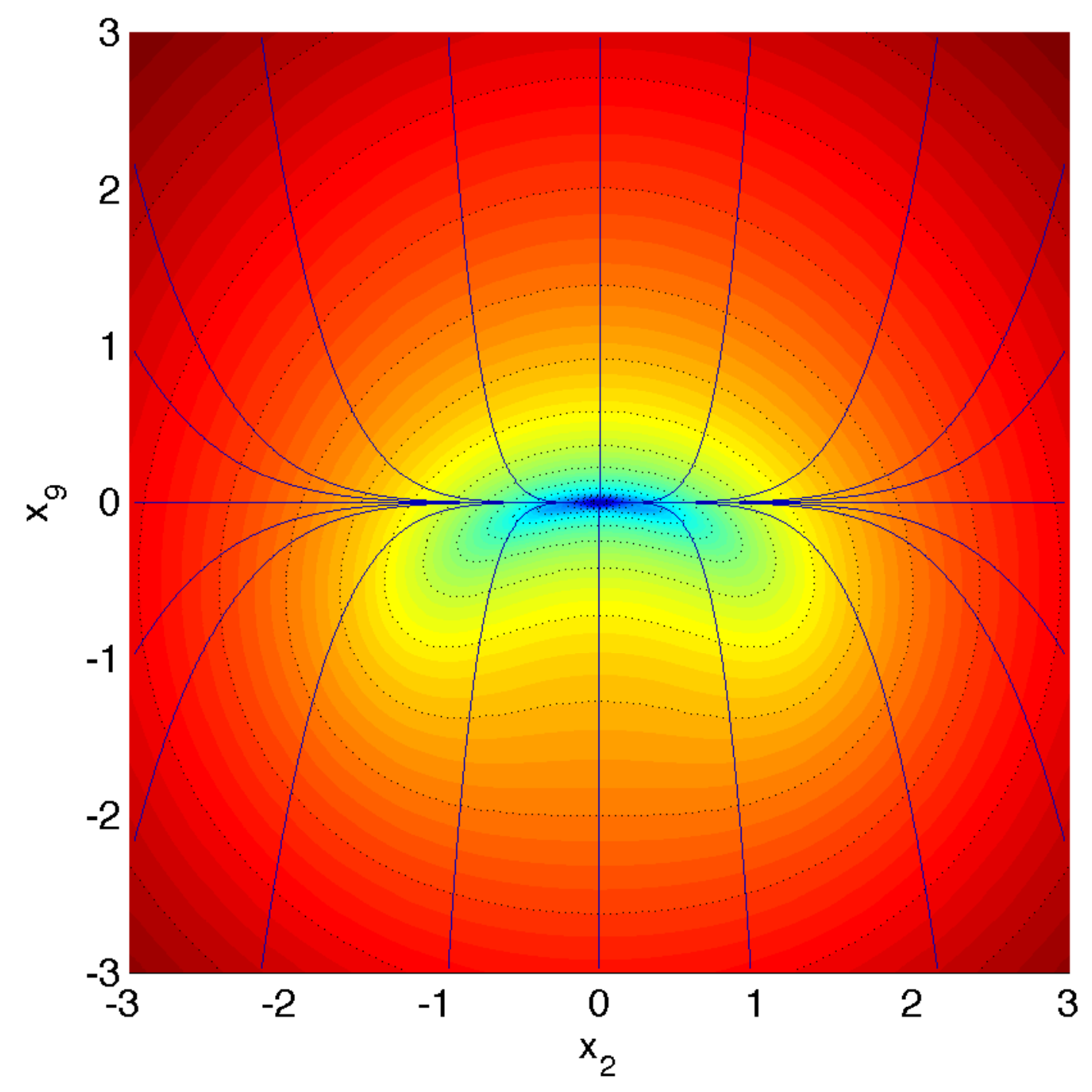}%
    }%
}~~%
  \caption[]{Selection of cross-sectional plots of the phase space for the nonlinear system $\dot a = L a + N(a)a$ with \mbox{$\Reynolds = 54$}.  Each figure shows a cut through the origin of the phase space for some $(a_i,a_j)$ pair.  Shading indicates the relative value of the Lyapunov function $V$ computed for Case 5 in Table \ref{tab:LyapResults}, with dashed black lines indicating level sets.  Streamlines of the nonlinear system, projected onto the cutting plane, are shown in solid blue.  
}%
\label{phasePlot}%
~\\[-1ex]
\hrulefill
\end{figure}

A selection of cross-sectional plots showing the phase space of the system $\dot a = f(a) = L(a) + N(a)a$ and level sets of the computed function $V$ for Case 5 are shown in Figure~\ref{phasePlot}.  The following features are of interest:
\begin{itemize} 
		\item Figures \ref{phasePlot}.\subref{phasePlot:13} and \ref{phasePlot}.\subref{phasePlot:19} show clearly that level sets of $V$ are centered on $a = 0$ for $\norm{a} \ll 1$, and are centered on $a = c = (1,\,0,\,\cdots,0)$ for $\norm{a} \gg 1$.
		\item Level sets of $V$ are radially symmetric in \ref{phasePlot}.\subref{phasePlot:67}.  Streamlines of $f$ approach tangents to these level sets as $\Reynolds \to \infty$.  Radial symmetry is also clear in \ref{phasePlot}.\subref{phasePlot:29} for $\norm{a} \gg 1$.
		\item There is significant warping of the level sets of $V$ (with respect to a radially symmetric energy function) in \ref{phasePlot}.\subref{phasePlot:19} and \ref{phasePlot}.\subref{phasePlot:29}.  Optimizing $V$ over a high-order polynomial allows increased freedom to shape these sets, thereby increasing the range of values $\Reynolds$ for which $f$ can be proven stable.
\end{itemize}

\subsection{Lyapunov Functions with Higher-Order Variable Terms}
\begin{figure}[!htpb]\centering
\scalebox{0.4}
{{\includegraphics[clip=true,trim = 25mm 5mm 25mm 3mm]{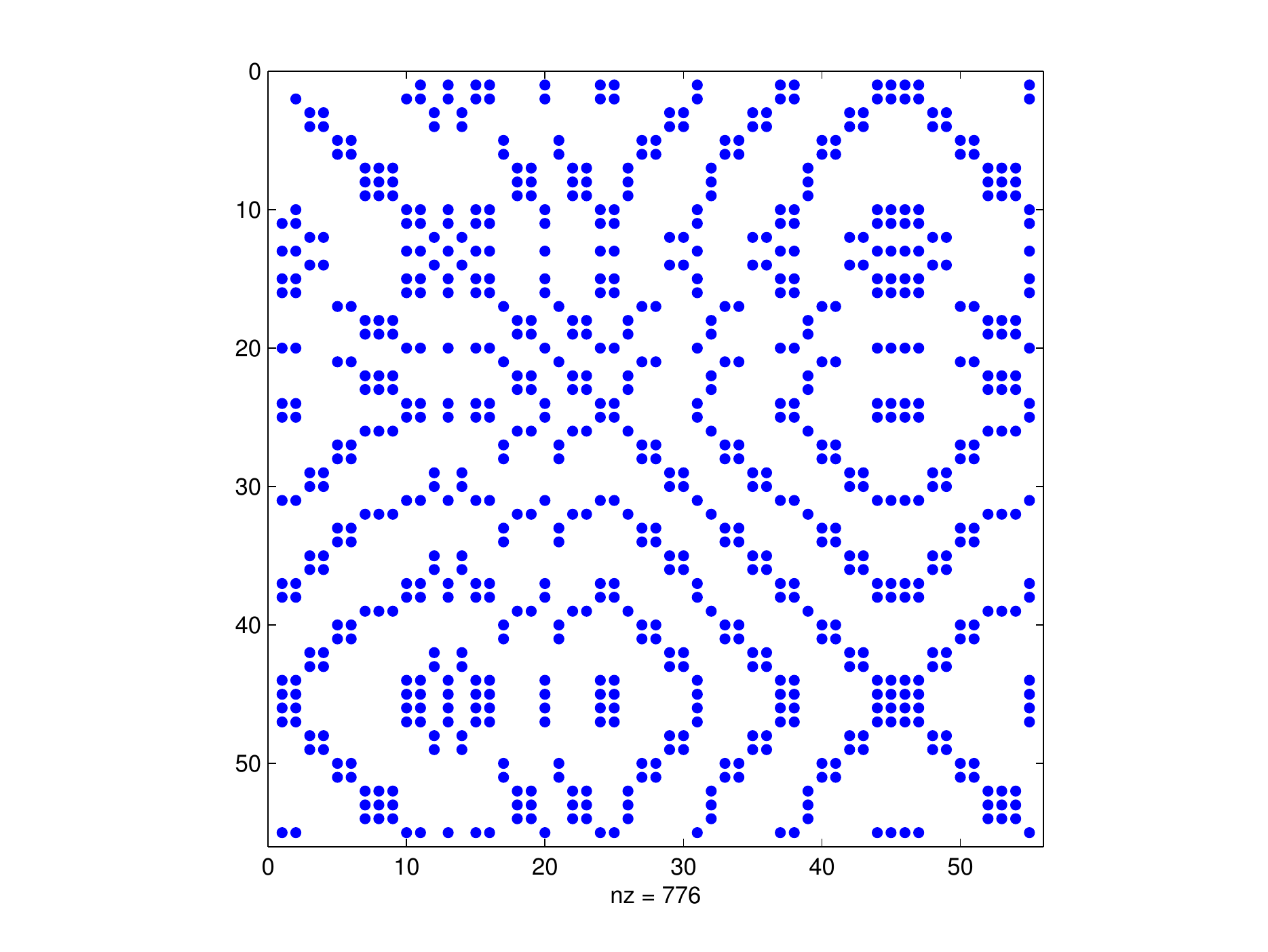}}}
{\raisebox{0.13\textheight}{$\xrightarrow[\text{reordering}]{\text{~~sparse~~}}$~~}}
\scalebox{0.4}
{{\includegraphics[clip=true,trim = 25mm 5mm 25mm 3mm]{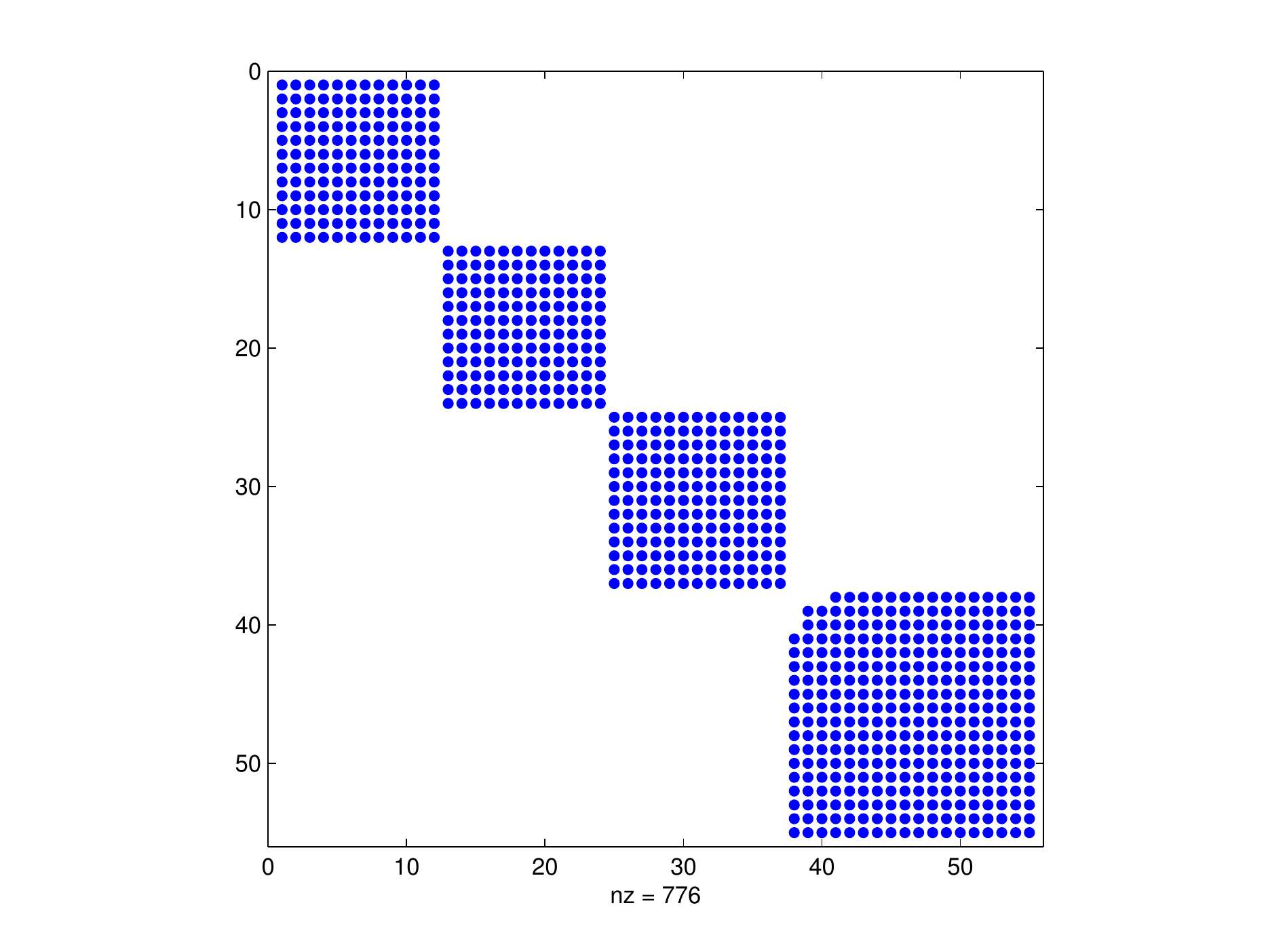}}}
\caption{Case 4: Reordering of monomial terms in $m_2(a)$ reveals block-diagonal structure in $P$.}\label{fig:ord4SparsityA}
\end{figure}

Direct solution of problem \eqref{SDP} for Lyapunov functions with variables terms of degree greater than four is more problematic given the long computation times required.  However, it is possible that close scrutiny of the results from cases 4--5 may give some indication of appropriate sparsity structures that may be exploited.  For example, as shown in Figure~\ref{fig:ord4SparsityA}, the sparsity pattern of the matrix $P$ featuring in the term $m_2(a)\tpose P m_2(a)$ in case 4 can be reordered to block diagonal form.  Note in particular that the upper-left hand corner of the unordered matrix in 
Figure~\ref{fig:ord4SparsityA}, which corresponds to the second-order terms in $m_2(a)\tpose P m_2(a)$, adopts an identical sparsity pattern to the solution of \eqref{ex:Psparsity}.
The resultant reordering is such that one can rewrite the variable component of $V$ as  
\[
m_2(a)\tpose P m_2(a) = 
\begin{bmatrix}
\tilde m^1_2(a)\\
\tilde m^2_2(a)\\
\tilde m^3_2(a)\\
\tilde m^4_2(a)
\end{bmatrix}\tpose
\begin{bmatrix}
\tilde P_1\\
& \tilde P_2\\
&& \tilde P_3\\
&&& \tilde P_4
\end{bmatrix}
\begin{bmatrix}
\tilde m^1_2(a)\\
\tilde m^2_2(a)\\
\tilde m^3_2(a)\\
\tilde m^4_2(a)\\
\end{bmatrix},
\] 
where each of the matrices $\tilde P_i$ is symmetric and dense and the reordered and partitioned monomial terms $\tilde m_2^i(a)$ are defined as 
\begin{equation}
\begin{aligned}
\tilde m_2^1(a) &\eqdef (a_2, a_3, a_1a_2, a_1a_3, a_4a_6, a_5a_6, a_4a_7, a_5a_7, a_4a_8, a_5a_8, a_2a_9, a_3a_9)\tpose\\
\tilde m_2^2(a)& \eqdef (a_4, a_5, a_1a_4, a_1a_5, a_2a_6, a_3a_6, a_2a_7, a_3a_7, a_2a_8, a_3a_8, a_4a_9, a_5a_9)\tpose \\
\tilde m_2^3(a)& \eqdef
 (a_6, a_7, a_8, a_2a_4, a_3a_4, a_2a_5, a_3a_5, a_1a_6, a_1a_7, a_1a_8, a_6a_9, a_7a_9, a_8a_9)\tpose \\
\tilde m_2^4(a)& \eqdef
(1, a_1, a_9, a_2a_3, a_4a_5,  a_6a_7,  a_6a_8, a_7a_8,  a_1a_9, a_1^2, a_2^2,  a_3^2, a_4^2, a_5^2, a_6^2, a_7^2, a_8^2, a_9^2 )\tpose.
\end{aligned}
\!\!\!\!\!\!\!
\end{equation}

\section{Conclusions}

A new method for analyzing the global stability of a fluid flows has been proposed. This method requires only the solutions of linear eigenvalue problems for systems of linear partial differential equations, combined with a nonlinear analysis of a system governed by nonlinear ordinary differential equations, which can be treated using the polynomial sum-of-squares approach. The method is proven always to yield results that at least as good as classical energy methods, in the sense that if the global stability of a particular flow can be proved by the energy method then it also can be proved by the new method. 

Moreover, if the flow remains globally stable for values of the Reynolds number in a range extending beyond the maximum for which global stability can be proved by the energy method, then a polynomial Lyapunov function is still guaranteed to exist in at least part of this extended range. The methods proposed in this paper can then apply, provided that the appropriate quantities can be expressed as sums-of-squares.

Application of our method to a finite dimensional model system suggests that using the proposed method might allow proving global stability for much higher values of the Reynolds number than the values for which this can be done by the energy method.  It would of course be of interest to extend such an example to the infinite dimensional case using the techniques of Section \ref{sec:InfiniteDim}.

\bibliographystyle{plain}
\bibliography{bibfile}    

\begin{appendix}

\section{Calculating the norm bounding coefficients $c_i.$}\label{app:NormBoundConstants}

In this appendix we will prove \eqref{lem:normSquareBounds}. The claim is that there exist constants $(c_1,c_2,c_3) \ge 0$ such that
\begin{equation}\label{app:eqn:normSquareBound}
\norm{\Theta_a(\mbf u_s) + \Theta_b(\mbf u_s,a) + \Theta_c(\mbf u_s)}^2
\le 
c_1 q^2 + c_2 \norm{a}^2q^2 + c_3 q^4 \,\,\forall a \in \Re^n, \forall \mbf{u}_s \in \mcf{S}.
\end{equation}
First note that for any norm, 
\[
\norm{u + v }^2 \le \norm{u + v}^2 + \norm{u-v}^2
= 2\left(\norm{u}^2 + \norm{v}^2\right).  
\]

We can apply this inequality to \eqref{app:eqn:normSquareBound} and then compute bounds for each of the terms in turn.

\subsection*{Part I -- Computing a bound on $\norm{\Theta_a(\mbf u_s)+\Theta_b
(\mbf u_s)}^2$:}

\[
\norm{\Theta_a(\mbf u_s)+\Theta_b(\mbf u_s)}^2=\sum_i \innerprod{\mbf u_s}{\mbf h_i+\mbf h_{ij}a_j}^2 \le I(a)\norm{\mbf u_s}^2
\]
where
\begin{equation}\label{eqn:Supremum}
I(a)=\sup_{\mbf u_s}\frac{\sum_i \innerprod{\mbf u_s}{\mbf h_i+\mbf h_{ij}a_j}^2}{\norm{\mbf u_s}^2}.
\end{equation}

Let $\tilde{\mbf h}_i =\mbf h_i+\nabla \phi_i+\beta_{ik}\mbf u_k$ be a projection of $h_i$ on the solenoidal subspace orthogonal to all $\mbf u_i$, so that $\nabla \tilde{\mbf h}_i=0$ and $\innerprod{\tilde{\mbf h}_i}{u_k}=0$. Then the scalar fields $\phi_i$ should satisfy the Poisson equation $\nabla^2\phi_i= -\nabla\cdot\mbf h_i.$ We will make the solution unique by imposing the boundary condition 
\[
\left.\frac{\partial \phi_i}{\partial n}\right|_{\partial\Omega}=-(\mbf h_i,\mbf n),
\]
where $\mbf n$ is a unit vector normal to the boundary $\partial\Omega.$ This condition ensures that the normal component of $\tilde{\mbf h}_i$ equals zero: $(\mbf h_i,\mbf n)=0.$ Similarly, we also introduce solenoidal fields $\tilde{\mbf h}_{ij}a.$

Strictly speaking, in \eqref{eqn:Supremum} the solenoidal velocity field $\mbf u_s$ should also satisfy the full boundary condition $\left.\mbf u_s\right|_{\partial\Omega}=0.$ However, this condition can be relaxed to $\left.(\mbf u_s,\mbf n)\right|_{\partial\Omega}=0$ without changing the value of $I(a),$ since for any vector field $\mbf v$ satisfying only $\left.(\mbf v,\mbf n)\right|_{\partial\Omega}=0$ one can easily construct a solenoidal field $\mbf v'$ satisfying $\left.\mbf v'\right|_{\partial\Omega}=0$ and such that $\innerprod{\mbf v-\mbf v'}{\mbf h_i}$ is arbitrarily small. With this relaxation one can represent $\mbf u_s$ in \eqref{eqn:Supremum} as
\[
\mbf u_s=U_i\tilde{\mbf h}_i+U_{ij}\tilde{\mbf h}_{ij}+\hat{\mbf u}_s,
\] 
where $\hat{\mbf u}_s$ is orthogonal to all $\tilde{\mbf h}_i$ and $\tilde{\mbf h}_{ij}.$ Obviously, $\hat{\mbf u}_s$ does not contribute to the numerator of \eqref{eqn:Supremum} but can only increase the denominator. Therefore, for calculating the supremum one can take $\hat{\mbf u}_s=0.$ Then \eqref{eqn:Supremum} becomes
\begin{equation}\label{eqn:SupremumFD}
I(a)=\sup_{U_k,U_{kl}}\frac{\sum_i \innerprod{U_k\tilde{\mbf h}_k+U_{kl}\tilde{\mbf h}_{kl}}{\tilde{\mbf h}_i+\tilde{\mbf h}_{ij}a_j}^2}{\norm{U_k\tilde{\mbf h}_k+U_{kl}\tilde{\mbf h}_{kl}}^2}.
\end{equation}

Since
\[
\sum_i \frac{\innerprod{U_k\tilde{\mbf h}_k+U_{kl}\tilde{\mbf h}_{kl}}{\tilde{\mbf h}_i+\tilde{\mbf h}_{ij}a_j}^2}{\norm{U_k\tilde{\mbf h}_k+U_{kl}\tilde{\mbf h}_{kl}}^2}
\le
\sum_i \frac{\innerprod{\tilde{\mbf h}_i+\tilde{\mbf h}_{ij}a_j}{\tilde{\mbf h}_i+\tilde{\mbf h}_{ij}a_j}^2}{\norm{\tilde{\mbf h}_i+\tilde{\mbf h}_{ij}a_j}^2}=\sum_i \norm{\tilde{\mbf h}_i+\tilde{\mbf h}_{ij}a_j}^2
\]
both the existence of the supremum and of the coefficients $c_1$ and $c_2$ such that $I(a)\le (c_1+c_2\norm{a}^2)/2$ becomes obvious. Note that the numerator of the estimate \eqref{eqn:SupremumFD} is a fourth-order polynomial in coefficients $U_k,$ $U_{kl},$ and $a_j,$ and that polynomial is quadratic separately in $U_k,$ $U_{kl},$ and $a_j,$ which implies that efficient numerical methods can be found for determining $c_1$ and $c_2.$    

\subsection*{Part II -- Computing a bound on $\norm{\Theta_c{(\mbf u_s)}}$:}

Calculation of this bounds is straightforward, and amounts to finding a solution to a generalized eigenvalue problem using a symmetric version of the operator $\nabla \mbf u_i$, in a manner similar to that used to calculate the bound $\lambda$ in~\eqref{eqn:EigenEnergy}.    This results in a constant $c_{3i}$ such that
\[
\norm{\Theta_{ci}(\mbf u_s)} \le c_{3i} q^2,
\]
so that $c_3 = 2\left(\sum\nolimits_i {\displaystyle c_{3i}^2}\right)$ provides a conservative bound.  
This completes the proof of \eqref{app:eqn:normSquareBound}.

\section{Existence of a polynomial Lyapunov functional for Reynolds numbers greater than the energy stability limit}\label{app:existence}

In general a flow can, and often will, remain globally stable in a certain range of Reynolds numbers greater than the maximum Reynolds number $\ReynoldsEnergyLim$ for which global stability can be proved by the energy method. We will prove now that at least in some part of this range there exists a Lyapunov functional that is polynomial in $a$ and $q,$ i.e.\ of a form suitable for the proposed method. To this end we will simply give an explicit expression for a function which is polynomial in $a$ and $q^2$ and which satisfies all the conditions for Lyapunov functionals. For this purpose we will assume that $\mbf u_i$ are chosen as the eigenfunctions of~(\ref{eqn:EigenEnergy}), and note that this is always possible. As a result, we will have $\chi=0$ as explained in Subsection~\ref{DescriptionasanUncertainSystem}. 
We will, of course, use explicitly also the assumption that the flow is indeed globally stable in some vicinity of the energy stability limit. More precisely, we require that 
\begin{equation}\label{eqn:Condition}
\sum_{i=2}^n  f_i^2/[(1+2E_0)E_0]
\end{equation}
 is strictly positive when $a_2=a_3=\dots=a_n=0$ for all values of $a_1$ including zero and infinity. If this were not true at some finite $a_1=a_{1s}$ for any $n$ that would mean that $a_{1s}\mbf u_1$ is a steady non-zero solution, which contradicts the assumption of the basic flow being globally stable%
\footnote{
Note that this argument is valid only if $n>1,$ and indeed one can check that for $n=1$ our approach will not work for $\Reynolds> \ReynoldsEnergyLim$ because the one-dimensional system turns out to be unstable.   It can be shown that in the general case $n$ has to be greater than the number of positive eigenvalues $\lambda_i.$}.
 This means that this condition can be satisfied by selecting large enough $n.$ To ensure the positiveness of (\ref{eqn:Condition}) at the origin we will additionally require that the flow is linearly asymptotically stable, which is not much of a loss of generality for the flow which is already assumed to be globally stable. To ensure the positiveness of \eqref{eqn:Condition} with $a_1$ at infinity, we will also assume that $\innerprod{\mbf u_i}{\mbf u_1\cdot \nabla\mbf u_1}\ne0$ for at least one $i$ in $2,\dots,n.$ If $\mbf u_1\cdot \nabla\mbf u_1\ne 0,$ which is often the case, this assumption can always be ensured to be true by selecting a sufficiently large $n.$ 
Finally, we will also assume that the eigenfunctions $\mbf u_i$ and eigenvalues $\lambda_i$  are continuous functions of $\Reynolds$  in some vicinity of the energy stability limit.  

The basic idea is simple. At $\Reynolds=\ReynoldsEnergyLim $ energy satisfies all the conditions for being a Lyapunov function everywhere except along $\mbf u_1,$ where its time derivative is zero. We add to it a small correction constructed in such a way that its time derivative is negative along $\mbf u_1.$ The modified expression will therefore be a Lyapunov function at $\Reynolds= \ReynoldsEnergyLim$ and, by continuity, also in its vicinity. In the sequel, we exploit this idea in a more formal proof. 

We propose the following expression for the Lyapunov function
\begin{equation}\label{eqn:Vmodified}
V=E_0+q^2+(E_0+q^2)^2-\varepsilon \sum_{i=2}^na_if_i,
\end{equation}
where $f(a)$ is defined by (\ref{eqn:sysNominal}) and $\epsilon$ is a small positive value.
Note that the sum in \eqref{eqn:Vmodified} starts at $i=2$. Note also that $f_i$ are of order $a$ for small $a$ and of order $a^2$ for large $a.$  

We first provide conditions under which this choice of $V$ is positive-definite.  
As required, $V=0$ if both $a=0$ and $q=0.$ Otherwise, $V>0$ provided $\varepsilon > 0$ and 
$$\varepsilon<\min \frac{E_0+E_0^2}{|\sum_{i=2}^na_if_i|}.$$ This ratio is easily shown to be bounded on an open ball of sufficiently small (large) radius.  The fraction is continuous outside (inside) of such a ball.  Existence of a positive minimizer is therefore ensured by the extreme value theorem.


%

We next develop conditions under which $\dot V$ is negative definite.  
Following the same arguments as in Subsection~\ref{ConversiontoaSumofSquaresproblem}, but noticing that here $\chi=0$ by assumption, we arrive at the following 
sufficient condition for satisfaction of the Lyapunov condition \eqref{eqn:lyapCondition}:
\begin{equation}\label{eqn:SuffCond} 
\left(\frac{\partial V}{\partial a} f(a) + 
\frac{\partial V}{\partial(q^2)} \cdot2\kappa_s q^2 \right)
< - 
\left|
\frac{\partial V}{\partial a} - \frac{\partial V}{\partial(q^2)}a\tpose 
\right|
p_1^\half(a,q),
\,\, \forall(a,q) \neq 0.
\end{equation}

From the definitions of $E_0$ and $\lambda_i$ and from the energy equation (\ref{eqn:EnergyEquation}) it follows that
$$
\frac{\partial E_0} {\partial a}\cdot f=\sum_{i=1}^n\lambda_ia_i^2.
$$
  
Using also the definition of $p_1$ in \eqref{lem:normSquareBounds}, the inequality \eqref{eqn:SuffCond} can be rewritten as 
\begin{multline} 
\frac{
\lambda_1a_1^2+2\kappa_s q^2+\sum_{i=2}^n\lambda_ia_i^2
}{E_0+q^2}+\qquad\qquad\qquad\qquad\qquad\qquad\qquad\qquad\qquad\qquad\qquad\qquad
\notag\\
\frac{
-\varepsilon\sum_{i=2}^n \left(f_i^2+a_i\sum_{j=1}^n \frac{\partial f_i}{\partial a_j}f_j\right) 
+\varepsilon\left|\frac{\partial\ }{\partial a}\sum_{i=2}^n a_if_i 
\right||q|(c_1 + c_2 \norm{a}^2 + c_3 q^2)^\half
}{(1+2E_0+2q^2)(E_0+q^2)}<0.
\label{eqn:SuffCondWithEpsilon}
\end{multline}

If we define
\begin{align}
\tilde a_i&\eqdef a_i/(E_0+q^2)^\half,\qquad
\tilde q \eqdef q/(E_0+q^2)^\half,\\ 
\tilde f_i &\eqdef f_i/[(1+2E_0+2q^2)(E_0+q^2)]^\half,\\
F_i &\eqdef F_i(a,q) \eqdef \frac{\sum_{j=1}^n \frac{\partial f_i}{\partial a_j}f_j}{{(1+2E_0+2q^2)(E_0+q^2)^\half}},\\
F_q &\eqdef F_q(a,q) \eqdef \left|\frac{\partial\ }{\partial a}\sum_{i=2}^n a_if_i 
\right|
\frac{(c_1 + c_2 \norm{a}^2 + c_3 q^2)^\half}{{(1+2E_0+2q^2)(E_0+q^2)^\half}},
\end{align}

then \eqref{eqn:SuffCondWithEpsilon} can be written more compactly as
\begin{equation}\label{eqn:SuffCondWithTildes}
\lambda_1\tilde a_1^2+ 2\kappa_s \tilde q^2+\sum_{i=2}^n\lambda_i\tilde a_i^2
-\varepsilon\sum_{i=2}^n \left(\tilde f_i^2+\tilde a_iF_i\right) 
+\varepsilon|q|F_q<0.
\end{equation}

Using the identity $xy=(Ax)^2-(Ax-y/(2A))^2+y^2/(4A^2)$, one can rewrite (\ref{eqn:SuffCondWithTildes}) as
\begin{align}
\lambda_1\tilde a_1^2+(2\kappa_s+\varepsilon A^2) \tilde q^2+
\sum_{i=2}^n(\lambda_i+\varepsilon A^2)\tilde a_i^2
-\varepsilon\sum_{i=2}^n \tilde f_i^2
+
\frac{\varepsilon}{4A^2}\left(F_q^2+\sum_{i=2}^n F_i^2\right) 
-
\notag
\\
\varepsilon\sum_{i=2}^n\left(A\tilde a_i+\frac{F_i}{2A}\right)^2
-
\varepsilon\left(A q-\frac{F_q}{2A}\right)^2<0.\label{eqn:SuffCondWithSquares}
\end{align}

The last two terms of \eqref{eqn:SuffCondWithSquares} are non-positive, so a sufficient condition to satisfy \eqref{eqn:SuffCondWithSquares} is
\begin{align}
D(\Reynolds)=\lambda_1\tilde a_1^2+(2\kappa_s+\varepsilon A^2+\varepsilon) \tilde q^2+
\sum_{i=2}^n(\lambda_i+\varepsilon A^2+\varepsilon)\tilde a_i^2
-\varepsilon\left[q^2+\sum_{i=2}^n (\tilde f_i^2+\tilde a_i^2)\right.
\notag\\
-\left.
\frac{1}{4A^2}\left(F_q^2+\sum_{i=2}^n F_i^2\right)\right]<0.
\label{eqn:SuffSuffCond}\end{align}

Note that if $\Reynolds= \ReynoldsEnergyLim$, then $\lambda_1=0,$ $\kappa_s<0$ and $\lambda_i<0$, with both inequalities strict. 
Satisfaction of \eqref{eqn:SuffSuffCond} is therefore ensured if both
$$
\varepsilon<-\min\{2\kappa_s,\lambda_i\}/(A^2+1)
$$ 
and $A$ is chosen sufficiently large to guarantee that
\[
q^2+\sum_{i=2}^n (\tilde f_i^2+\tilde a_i^2)-
\frac{1}{4A^2}\left(F_q^2+\sum_{i=2}^n F_i^2\right) > 0
\]
for all $(a,q) \neq 0$.  
A suitable choice is 
\begin{equation}\label{eqn:A2boundfrac}
A^2>\frac14\sup\frac{ F_q^2+\sum_{i=2}^n F_i^2}{q^2+\sum_{i=2}^n (\tilde f_i^2+\tilde a_i^2)},
\end{equation}
provided that this fraction can be shown to be bounded for all $(a,q)$.
 
We first consider whether the numerator of \eqref{eqn:A2boundfrac} is bounded above. The functions $F_i$ and $F_q$ are easily shown to be bounded inside (outside) an open ball of sufficiently small (large) radius, and are continuous elsewhere.  Boundedness over all $(a,q)$ then follows from the boundedness theorem.   

The denominator in \eqref{eqn:A2boundfrac} is bounded below by a strictly positive value, because we already assumed that (\ref{eqn:Condition}) is bounded below by a strictly positive value when $a_2=a_3=\dots=a_n=0.$

The above argument ensures that $D(\ReynoldsEnergyLim)<0.$
Then from our continuity assumptions it follows that $D(\Reynolds)<0$ in at least some vicinity where $\Reynolds> \ReynoldsEnergyLim,$ thus completing the proof.





\section{System Dynamics for Numerical Example in \S\ref{sec:example}}\label{app:example}

The shear flow model used in the example in Section \ref{sec:example} is taken directly from~\mbox{\cite{MFE04,MFE05}}.   We include here the basis functions $\mbf{u}_i$ and resulting ODE system from~\cite{MFE04,MFE05} for easy reference.

The model uses the following 9--dimensional basis of mutually orthogonal, solenoidal basis functions:
\begin{alignat*}{2}
\mbf{u}_1 &= \begin{bmatrix}
				\sqrt{2} \sin (\pi y/2) \\
				0 \\
				0
			\end{bmatrix}\!,\,
			&
\mbf{u}_2 &= 
			 \begin{bmatrix}
			 	\cos^2 (\pi y/2)\cos(\gamma z) \\
			 	0 \\
			 	0
			 \end{bmatrix}\cdot 
			 \frac{4}{\sqrt{3}}, \\
\mbf{u}_3 &= 
			 \begin{bmatrix}
				0\\
				2\gamma\cos(\pi y/2) \cos(\gamma z)\\
				\hphantom{1}\pi  \sin(\pi y/2) \sin(\gamma z)
			 \end{bmatrix}\cdot 
			 \frac{2}{\sqrt{4\gamma^2+\pi^2}}
			 &
\mbf{u}_4 &= 
			 \begin{bmatrix}
			 	0\\ 
				0\\ 
				\cos (\alpha  x) \cos^2(\pi y/2)\\
			 \end{bmatrix}\cdot 
			 \frac{4}{\sqrt{3}}, \\
\mbf{u}_5 &= \begin{bmatrix}
				0\\
				0\\
				2 \sin(\alpha  x) \sin(\pi y/2)
			 \end{bmatrix}\\
\mbf{u}_6 &= 
			 \begin{bmatrix}
			 	-\gamma  \cos (\alpha  x) \cos^2(\pi y/2) \sin(\gamma z)\\
				0\\
				\alpha   \sin(\alpha  x) \cos^2(\pi y/2) \cos (\gamma z)\\
			  \end{bmatrix} \cdot 
			  \hsmash{\frac{4\sqrt{2}}{\sqrt{3(\alpha  ^2 + \gamma ^2)}},}\\
\mbf{u}_7 &= 
			 \begin{bmatrix}
			 	\hphantom{-}\gamma  \sin(\alpha  x) \sin(\pi y/2) \sin(\gamma z)\hphantom{{,}^2}\\ 
				0\\
				\alpha   \cos (\alpha  x) \sin(\pi y/2) \cos (\gamma z)
			  \end{bmatrix}\cdot 
			  \hsmash{\frac{2\sqrt{2}}{\sqrt{(\alpha  ^2 + \gamma ^2)}},}\\
\mbf{u}_8 &= 
			  \begin{bmatrix}
			  \pi \alpha   \sin(\alpha  x) \sin(\pi y/2) \sin(\gamma z)\\
			  2(\alpha ^2 + \gamma^2) \cos (\alpha  x) \cos (\pi y/2) \sin(\gamma z) \\
			   -\pi \gamma  \cos (\alpha  x) \sin(\pi y/2) \cos (\gamma z)
			  \end{bmatrix}\cdot 
			 N_8, & \quad 
			 \mbf{u}_9 &= \begin{bmatrix}
			  \sqrt{2} \sin(3\pi y/2)\\ 
			  0\\
			  0
\end{bmatrix},
\end{alignat*}
where 
$\alpha = 2\pi/L_x$, $\beta = \pi / 2$, $\gamma = 2\pi /L_z$ and
\[
N_8 ={\frac{2\sqrt{2}}{\sqrt{(\alpha^2 + \gamma^2)(4\alpha^2 + 4\gamma^2 + \pi^2) }}.}
\]   It is easily verified that $\mbf{w} = \mbf{u}_1$ is a laminar solution of the Navier-Stokes equation \eqref{eqn:NS} when the volume force is
\[
\mbf{f} = \Bigl(\frac{\sqrt{2}\pi^2}{4\Reynolds}\sin(\pi y/2);0;0\Bigr). 
\]

\newcommand{\kab}{\kappa_{\alpha\beta}}
\newcommand{\kbg}{\kappa_{\beta\gamma}}
\newcommand{\kag}{\kappa_{\alpha\gamma}}
\newcommand{\kabg}{\kappa_{\alpha\beta\gamma}}

These basis functions can then be expanded via Galerkin projection into a nonlinear system of ODEs as described in Section \ref{sec:finiteDimProblems}.  Define the following notation for neatness:
\begin{align*}
\kab \eqdef \sqrt{\alpha^2 + \beta^2}, \,\,
\kag \eqdef \sqrt{\alpha^2 + \gamma^2}, \,\,
\kbg \eqdef \sqrt{\beta^2 + \gamma^2}, \,\,\textrm{and }
\kabg \eqdef \sqrt{\alpha^2 + \beta^2+ \gamma^2}.
\end{align*}
Then the ODE in the numerical example is the form \eqref{eqn:sysNominal:alt}, with
$c = \operatorname{diag}(1,\,\,0,\dots, 0)$,
{
\small
\[
\Lambda = -\operatorname{diag}\left(%
\beta^2,\,\,										
\frac{4\beta^2}{3} + \gamma^2,\,\,	 			
\kbg^2,\,\,										
\frac{3 \alpha^2 + 4 \beta^2}{3},\,\,				
\kab^2,\,\,										
\frac{3 \alpha^2 + 4 \beta^2 + 3 \gamma^2}{3},\,\,	
\kabg^2,\,\,				
\kabg^2,\,\,				
9\beta^2
\right),
\]
} 
and 
{
\small
\begin{align*}
[N(a)a]_1&=
	\sqrt{\frac{3}{2}}\frac{\beta\gamma}{\kbg} a_2 a_3 - 
	\sqrt{\frac{3}{2}} \frac{\beta\gamma}{\kabg} a_6 a_8 
\\
[N(a)a]_2&=
	\frac{10}{3\sqrt{6}}\frac{\gamma^2}{\kag}a_4 a_6 - 
	\frac{\gamma^2}{\sqrt{6} \kag} a_5 a_7-
	\frac{\alpha\beta\gamma}{\sqrt{6}\kag \kabg} a_5 a_8 - 
	\sqrt{\frac{3}{2}} \frac{\beta\gamma}{\kbg} \left(a_1 a_3 + a_3 a_9\right) 
\\
[N(a)a]_3&=
	\sqrt\frac{2}{3} \frac{\alpha\beta\gamma}{\kag \kbg}\left(a_5 a_6 + a_4 a_7\right) + 
	\frac{\beta^2 (3 \alpha^2 + \gamma^2) -3 \gamma^2 \kag^2}{\sqrt{6} \kag \kbg \kabg}a_4 a_8
\\
[N(a)a]_4&=
	 -\frac{\alpha}{\sqrt{6}} (a_1 a_5+a_5 a_9) - 
	 \frac{10}{3\sqrt{6}} \frac{\alpha^2}{\kag} a_2 a_6 - 
	 \sqrt{\frac{3}{2}} \frac{\alpha\beta\gamma}{\kag \kbg} a_3 a_7 - 
	 \sqrt{\frac{3}{2}} \frac{\alpha^2 \beta^2}{\kag \kbg \kabg} a_3 a_8
\\
[N(a)a]_5&=
	\frac{\alpha}{\sqrt{6}} (a_1 a_4 +a_4 a_9) +
	\sqrt{\frac{2}{3}} \frac{\alpha\beta\gamma}{\kag \kbg}a_3 a_6 +
	\frac{\alpha^2}{\sqrt{6} \kag}a_2 a_7 - 
	\frac{\alpha\beta\gamma }{\sqrt{6} \kag \kabg} a_2 a_8
\\
[N(a)a]_6&=
	\frac{10}{3\sqrt{6}}\frac{\alpha^2 - \gamma^2}{\kag}a_2 a_4 - 
	\sqrt{\frac{2}{3}}\frac{ 2\alpha\beta\gamma }{\kag \kbg}a_3 a_5 + 
	\frac{\alpha}{\sqrt{6}} (a_1 a_7 +a_7 a_9) + 
	\sqrt{\frac{3}{2}} \frac{\beta\gamma}{\kabg} (a_1 a_8 + a_8 a_9 )
\\
[N(a)a]_7&=
	\frac{\alpha\beta\gamma}{\sqrt{6} \kag \kbg} a_3 a_4 + 
	\frac{(-\alpha^2 + \gamma^2) }{\sqrt{6} \kag}a_2 a_5 - 
	\frac{\alpha}{\sqrt{6}} (a_1 a_6 + a_6 a_9)  
\\
[N(a)a]_8&=
	\frac{\gamma^2 (3 \alpha^2 - \beta^2 + 3 \gamma^2)}{\sqrt{6} \kag \kbg\kabg} a_3 a_4 + 
	\sqrt{\frac{2}{3}}\frac{ \alpha\beta\gamma }{\kag \kabg} a_2 a_5 
\\
[N(a)a]_9&=
	\sqrt{\frac{3}{2}} \frac{\beta\gamma}{\kbg} a_2 a_3 - 
	\sqrt{\frac{3}{2}} \frac{\beta\gamma}{\kabg} a_6 a_8,	
\end{align*}%
}
where $[N(a)a]_i$ is the $i^\text{th}$ component of $N(a)a$ and $\operatorname{diag}(\cdot)$ forms a diagonal matrix from its arguments.  
%
%

%
%
%
%
%
%
%

\end{appendix}
\end{document}